\documentclass[11pt]{article}
\usepackage{amsmath}
\usepackage{amstext,amsbsy}
\usepackage{epsfig,multicol}
\usepackage{amsthm}
\usepackage{amssymb,latexsym}
\usepackage{color}
\usepackage{eepic}

\newtheorem{thm}{Theorem}[section]

\newtheorem{exa}[thm]{Example}
\newtheorem{lem}[thm]{Lemma}

\newtheorem{prop}[thm]{Proposition}

\def\pf{\noindent{\it Proof.} }

\def\qed{\nopagebreak\hfill{\rule{4pt}{7pt}}\medbreak}
\def \se {\mathrm{se}}
\def \ne {\mathrm{ne}}

\def \inv {\mathrm{inv}}
\def \asc {\mathrm{coinv}}
\def \F {\mathbf{F}}
\def \s {\mathbf{s}}
\def \bve {\mathbf{e}}
\def \ve  {\varepsilon}

\newcommand {\auc} {\mathrm{auc}}
\newcommand {\buc} {\mathrm{buc}}
\newcommand {\ce} {\mathrm{ce}}
\def \M {\mathcal{M}}
\def \N {\mathcal{N}}
\def \R {\mathcal{R}}
\def \B {\mathcal{B}}

\begin{document}
\title{Mixed Statistics on $01$-Fillings of Moon Polyominoes}
\author{
William Y. C. Chen$^{1}$, Andrew Y. Z. Wang$^1$, Catherine H.
Yan$^2$ and \\
Alina F. Y. Zhao$^{1}$\\ [6pt]
$^{1}$Center for Combinatorics, LPMC-TJKLC\\
Nankai University, Tianjin 300071, P. R. China\\[5pt]
$^{2}$Department of Mathematics\\
Texas A\&M University, College Station, TX 77843, USA\\[5pt]
}
\date{}
\maketitle
\begin{abstract}

 We establish a stronger symmetry
between the numbers of northeast and southeast chains in the context
of 01-fillings of moon polyominoes.
 Let $\M$ be a moon polyomino with $n$ rows and
$m$ columns. Consider all the 01-fillings of $\M$ in which every row
has at most one $1$. We introduce four mixed statistics with respect
to a bipartition of rows or columns of $\M$. More precisely, let $S
\subseteq \{1,2,\ldots, n\}$ and $\mathcal{R}(S)$ be the union of
rows whose indices are in $S$. For any filling $M$, the top-mixed
(resp. bottom-mixed) statistic $\alpha(S; M)$ (resp. $\beta(S; M)$)
is the sum of the number of northeast chains whose top (resp.
bottom) cell is in $\mathcal{R}(S)$, together with the number of
southeast chains whose top (resp. bottom) cell is in the complement
of $\mathcal{R}(S)$. Similarly, we define the left-mixed and
right-mixed  statistics $\gamma(T; M)$  and $\delta(T; M)$, where
$T$ is a subset of  the column index set $\{1,2,\ldots, m\}$. Let
$\lambda(A; M)$ be any of these four statistics $\alpha(S; M)$,
$\beta(S; M)$, $\gamma(T; M)$ and $\delta(T; M)$, we show that the
joint distribution of the pair $(\lambda(A; M), \lambda(\bar A; M))$
is symmetric and independent of the subsets $S, T$. In particular,
the pair of statistics $(\lambda(A;M), \lambda(\bar A; M))$ is
equidistributed with $(\se(M),\ne(M))$, where $\se(M)$ and $\ne(M)$
are the numbers of southeast chains  and northeast chains of $M$,
respectively.
\end{abstract}



{\renewcommand{\thefootnote}{}
\footnote{\emph{E-mail addresses}: chen@nankai.edu.cn (W.Y.C.~Chen),
yezhouwang@mail.nankai.edu.cn (A.Y.Z.~Wang), \\
cyan@math.tamu.edu (C.H.~Yan), zfeiyan@mail.nankai.edu.cn
(A.F.Y.~Zhao). } }

\footnotetext[1]{The first, the second, and the fourth authors were
supported by the 973 Project, the PCSIRT Project of the Ministry of
Education, and the National Science Foundation of China.}

\footnotetext[2]{The third author was supported in part by NSF grant
\#DMS-0653846.}

\section{Introduction} \label{S:introduction}

Recently it is observed that the numbers of crossings and nestings
have a symmetric distribution over many families of combinatorial
objects,   such as  matchings and set partitions. Recall that a
matching of $[2n]=\{1, 2, \ldots, 2n\}$ is a partition of the set
$[2n]$ with the property that each block has exactly two elements.
It can be represented as a graph with vertices $1,2, \ldots, 2n$
drawn on a horizontal line in increasing order, where two vertices
$i$ and $j$ are connected by an edge if and only if $\{i, j\}$ is a
block. We say that two edges $(i_1, j_1)$ and $(i_2, j_2)$ form a
\emph{crossing} if $i_1 < i_2< j_1 < j_2$; they form a
\emph{nesting} if $i_1 < i_2 < j_2 < j_1$. The symmetry of the joint
distribution of crossings and nestings follows from the bijections
of de Sainte-Catherine, who also found the generating functions for
the number of crossings and the number of nestings. Klazar
\cite{klazar} further studied the distribution of crossings and
nestings over the set of matchings obtained from a given matching by
successfully adding  edges.

The symmetry between crossings and nestings was extended by Kasraoui
and Zeng \cite{Kasra06} to set partitions, and by Chen, Wu and Yan
\cite{Chen08b} to linked set partitions. Poznanovi\'{c} and Yan
\cite{PY09}  determined the distribution of crossings and nestings
over the set of partitions which are identical to a given partition
$\pi$ when restricted to the last $n$ elements.

Many classical results on enumerative combinatorics can be put in
the larger context of counting submatrices in fillings of certain
polyominoes. For example, words and permutations can be represented
as 01-fillings of rectangular boards, and general graphs can be
represented as $\mathbb{N}$-fillings of arbitrary Ferrers shapes,
which were studied by Kratthenthaler \cite{Krattenthaler} and de
Mier \cite{deMier06,deMier}. Other extensions include stack
polyominoes \cite{Jon05}, and moon polyominoes \cite{Rubey,
Kasra08}. In particular, crossings and nestings in matchings and set
partitions correspond to northeast chains and southeast chains of
length $2$ in a filling of polyominoes. The symmetry between
crossings and nestings  has been extended by Kasraoui \cite{Kasra08}
to 01-fillings of moon polyominoes where either every row has at
most one $1$, or every column has at most one $1$. In both cases,
the joint distribution of the numbers of northeast and southeast
chains can be expressed  as a product of $p,q$-Gaussian
coefficients. Other known statistics  on fillings of moon
polyominoes are the length of the longest northeast/southeast chains
\cite{Chen07,Krattenthaler,Rubey}, and the major index
\cite{Chen09}.

The main objective of this paper is to present a stronger symmetry
between the numbers of northeast and southeast chains in the context
of 01-fillings of moon polyominoes. Given a bipartition of the rows
(or columns) of a moon polyomino, we define four statistics by
considering mixed sets  of northeast and southeast chains according
to the bipartition.  Let $M$ be a 01-filling of a moon polyomino
$\M$ with $n$ rows and $m$ columns. These statistics are the
top-mixed and the bottom-mixed statistics $\alpha(S;M), \beta(S; M)$
with respect to a row-bipartition $(S, \bar S)$, and the left-mixed
and the right-mixed statistics $\gamma(T; M), \delta(T;M)$ with
respect to a column-bipartition $(T, \bar T)$. We show that
 for any of these four  statistics $\lambda(A; M)$, namely,
$ \alpha(S; M), \beta(S;M)$ for $S \subseteq [n]$ and $\gamma(T; M),
\delta(T; M)$ for $T \subseteq [m]$,  the joint distribution of the
pair $(\lambda(A; M), \lambda(\bar A; M))$ is symmetric and
independent of the subsets $S, T$. Consequently, we have the
equidistribution
\[
 \sum_{M} p^{\lambda(A;M)} q^{\lambda(\bar A; M)} = \sum_M p^{\se(M)}q^{\ne(M)},
\]
where $M$ ranges over all 01-fillings of $\M$ with the property that
either every row has at most one $1$, or every column has at most
one $1$, and $\se(M)$ and $\ne(M)$ are the numbers of southeast and
northeast chains of $M$, respectively.

The paper is organized as follows. Section 2  contains necessary
notation and the statements of the main results. In Section 3, we
explain how our results specialize to classical combinatorial
objects, including permutations, words, matchings, and set
partitions. We present the proofs of the main theorems in Section 4.
In Section 5, we show by bijections that these new statistics are
invariant under a permutation of columns or rows on moon
polyominoes.

\section{Notation and  the Main Results}

A \textit{polyomino} is a finite subset of $\mathbb{Z}^2$, where
every element of $\mathbb{Z}^2$ is represented by a square cell. The
polyomino is \textit{convex} if its intersection with any column or row
is connected.
 It is \textit{intersection-free} if every
two columns are comparable, i.e., the row-coordinates of one
column form a subset of those of the other column. Equivalently, it is
intersection-free if every two rows are comparable. A \textit{moon
polyomino} is a convex and intersection-free polyomino.

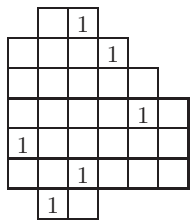
\begin{figure}[ht]
\begin{center}
\setlength{\unitlength}{0.4mm}
\begin{picture}(60,70)

\qbezier[500](10,0)(20,0)(30,0)   \qbezier[500](0,10)(30,10)(60,10)
\qbezier[500](0,20)(30,20)(60,20) \qbezier[500](0,30)(30,30)(60,30)
\qbezier[500](0,40)(30,40)(60,40) \qbezier[500](0,50)(25,50)(50,50)
\qbezier[500](0,60)(20,60)(40,60) \qbezier[500](10,70)(20,70)(30,70)
\qbezier[500](0,10)(0,35)(0,60)   \qbezier[500](10,0)(10,35)(10,70)
\qbezier[500](20,0)(20,35)(20,70) \qbezier[500](30,0)(30,35)(30,70)
\qbezier[500](40,10)(40,35)(40,60)\qbezier[500](50,10)(50,30)(50,50)
\qbezier[500](60,10)(60,25)(60,40)

\put(3,22){\footnotesize{1}} \put(13,2){\footnotesize{1}}
\put(23,12){\footnotesize{1}} \put(23,62){\footnotesize{1}}
\put(33,52){\footnotesize{1}} \put(43,32){\footnotesize{1}}
\end{picture}
\end{center}

\caption{A filling $M$  with $\mathbf{e}=(1,1,0,1,1,1,1)$ and
$\mathbf{s}=(1,1,2,1,1,0)$.}\label{polyomino}

\end{figure}

Given a moon polyomino $\mathcal {M}$, we assign  $0$ or $1$ to each
cell of $\mathcal {M}$ so that there is at most one $1$ in each row.
Throughout this paper we will simply use the term \emph{filling} to
denote such $01$-fillings. We say that a cell is empty if it is
assigned  $0$, and it is a $1$-cell otherwise. Assume $\mathcal {M}$
has $n$ rows and $m$ columns. We label the rows $R_1, \dots, R_n$
from top to bottom, and  the columns $C_1,\ldots,C_m$ from left to
right. Let $\bve=(\ve_1, \dots, \ve_n) \in \{0,1\}^n$ and
$\s=(s_1,\ldots,s_m) \in \mathbb{N}^m$ with \[ \sum_{i=1}^n \ve_i =
\sum_{j=1}^m s_j.\] We denote by $\mathbf{F}(\mathcal {M},\bve,
\mathbf{s})$ the set of fillings $M$ of $\mathcal {M}$ such that the
row $R_i$ has exactly $\ve_i$ many $1$'s, and the column $C_j$ has
exactly $s_j$ many $1$'s, for $ 1 \leq i \leq n$ and  $1\leq j\leq
m$. See Figure \ref{polyomino}  for an illustration.

A \emph{northeast} (resp. \emph{southeast}) \emph{chain} in a
filling $M$ of $\mathcal{M}$ is a set of two $1$-cells such that one
of them is strictly above (resp. below) and to the right of the
other and the smallest rectangle containing them is contained in
$\mathcal{M}$. Northeast (resp. southeast) chains will be called NE
(resp. SE) chains. The number of NE (resp. SE) chains of $M$ is
denoted by $\ne(M)$ (resp. $\se(M)$).

Let $\mathcal{R}$ be the set of rows of the moon polyomino $\M$. For
$S\subseteq [n]$, let \[ \mathcal{R}(S)= \bigcup_{i\in S} R_i.\]  We
say a $1$-cell is an $S$-\textit{cell} if it lies in
$\mathcal{R}(S)$. An NE chain is called a \emph{top  $S$-NE chain}
if its northeast $1$-cell is an $S$-cell. Similarly, an SE chain is
called a \emph{top  $S$-SE chain} if its northwest $1$-cell is an
$S$-cell. In other words, an NE/SE chain is a  top $S$-NE/SE chain
if the upper  $1$-cell of the chain is in  $\mathcal{R}(S)$.
Similarly, an NE/SE chain is a \emph{bottom $S$-NE/SE chain} if the
lower  $1$-cell of the chain is in  $\mathcal{R}(S)$.

Let $\bar{S}=[n]\setminus S$ be the complement of $S$.
Given a filling $M\in\mathbf{F}(\mathcal{M},\bve, \s)$,
we define the \emph{top-mixed  statistic $\alpha(S;M)$} and the \emph{bottom-mixed
statistic $\beta(S;M)$ with respect to $S$}  as
\begin{eqnarray*}
\alpha(S;M)&=&\#\{ \text{top } S\mbox{-NE chain of}~M\}+\#\{\text{top } \bar{S}
\mbox{-SE chain of}~M\}, \\[3pt]
\beta(S; M)&=& \#\{ \text{bottom } S\mbox{-NE chain of}~M\}+\#\{\text{bottom } \bar{S}
\mbox{-SE chain of}~M\}.
\end{eqnarray*}
See Example \ref{ex1} for some of these statistics on the filling $M$ in Figure \ref{polyomino}.

Let $F^t_S(p,q)$ and $F^b_S(p,q)$ be the bi-variate generating
functions for the pairs $(\alpha(S; M), \alpha(\bar S; M))$ and
$(\beta(S;M), \beta(\bar S; M))$ respectively, namely,
\[
 F^t_S(p,q)=\sum\limits_{M\in\mathbf{F}(\mathcal {M}, \bve, \s)}
p^{\alpha(S;M)} q^{\alpha(\bar{S};M)}
\quad \text{ and } \quad
F^b_S(p,q)=\sum\limits_{M\in\mathbf{F}(\mathcal {M}, \bve, \s)}
p^{\beta(S;M)} q^{\beta(\bar{S};M)}.
\]
Note that \[(\alpha(\emptyset; M), \alpha([n];
M))=(\beta(\emptyset;M), \beta([n];M))= (\se(M), \ne(M)).\] Our
first result is the following property.

\begin{thm}\label{thm-row}
$F^t_S(p,q)=F^t_{S'}(p,q)$ for any two subsets $S,S'$ of $[n]$. In
other words, the bi-variate generating function $F^t_S(p,q)$  does
not depend on $S$. Consequently,
\[F^t_S(p,q)=F^t_{\emptyset}(p,q)=\sum\limits_{M\in\mathbf{F}
(\mathcal {M},\bve, \s)}p^{\se(M)}q^{\ne(M)}
\]
is a symmetric function. The same statement holds for $F^b_S(p,q)$.
\end{thm}

We can also define the mixed statistics with respect to a subset of
columns. Let $\mathcal{C}$ be the set of columns of $\M$. For $T
\subseteq [m]$, let
\[ \mathcal{C}(T)= \bigcup _{j \in T} C_j .\]  An NE chain is
called  a \emph{left  $T$-NE chain}  if the southwest  $1$-cell of
the chain lies in  $\mathcal{C}(T)$. Similarly, an SE chain is
called a \emph{left  $T$-SE chain} if the northwest $1$-cell of the
chain lies in  $\mathcal{C}(T)$. In other words, an NE/SE chain is a
left $T$-NE/SE chain if its left $1$-cell is in $\mathcal{C}(T)$.
Similarly, an NE/SE chain is a \emph{right $T$-NE/SE chain} if its
right $1$-cell is in $\mathcal{C}(T)$.

Let $\bar{T} = [m]\setminus T$ be the complement of $T$.
For any filling $M$ of $\F(\M, \bve, \s)$,
we define the \emph{left-mixed  statistic $\gamma(T;M)$} and the \emph{right-mixed statistic
$\delta(T;M)$ with respect to $T$}  as
\begin{eqnarray*}
 \gamma(T; M)& = & \#\{\text{left }  T\mbox{-NE chain of}~M\}+\#\{\text{left } \bar{T}
\mbox{-SE chain of}~M\} , \\[3pt]
\delta(T; M) & =& \#\{\text{right }  T\mbox{-NE chain of}~M\}+\#\{\text{right } \bar{T}
\mbox{-SE chain of}~M\}.
\end{eqnarray*}

\begin{exa} \label{ex1}
\emph{Let $M$ be the filling in Figure \ref{polyomino}, where
$\ne(M)=6$ and $\se(M)=1$. Let $S=\{2,4\}$, i.e., $\R(S)$
contains the second and
the fourth rows. Then }
\[
 \alpha(S;M)=5, \quad \alpha(\bar S;M)=2, \quad
 \beta(S;M)=1,\quad \beta(\bar S;M)=6.
\]
\emph{Let $T=\{1, 3\}$, i.e., $\mathcal{C}(T)$ contains the first and the third columns.  Then }
\begin{equation*}
 \gamma(T;M)=4, \quad \gamma(\bar T; M)=3, \quad
\delta(T;M)=2, \quad \delta(\bar T;M)=5. \tag*{\rule{4pt}{7pt}}
\end{equation*}
\end{exa}

Let $G^l_T(p,q)$ and $G^r_T(p,q)$ be the bi-variate generating
functions of the pairs $(\gamma(T; M), \gamma(\bar{T}; M))$ and
$(\delta(T;M), \delta(\bar T; M))$ respectively, namely,
\[
G^l_T(p,q)=\sum\limits_{M\in\mathbf{F}(\mathcal {M},\bve, \s)}
p^{\gamma(T;M)} q^{\gamma(\bar{T};M)}
\quad \text{ and }  \quad
G^r_T(p,q)=\sum\limits_{M\in\mathbf{F}(\mathcal {M},\bve, \s)}
p^{\delta(T;M)} q^{\delta(\bar{T};M)}.
\]
Again note that
\[ (\gamma(\emptyset; M), \gamma([m]; M))=
(\delta(\emptyset;M), \delta([m]; M))=  (\se(M), \ne(M)).\] Our
second result shows that the generating function $G^l_T(p,q)$
possesses a similar property as $F^t_S(p,q)$.

\begin{thm}\label{thm-column}
$G^l_T(p,q)=G^l_{T'}(p,q)$ for any two subsets $T, T'$ of $[m]$. In
other words, the bi-variate generating function $G^l_T(p,q)$ does
not depend on $T$. Consequently,
\[
G^l_T(p,q)=G^l_{\emptyset}(p,q)=\sum\limits_{M\in\mathbf{F}
(\mathcal {M},\bve, \s)}p^{\se(M)}q^{\ne(M)}
\]
is a symmetric function. The same statement holds for $G^r_T(p,q)$.
\end{thm}

 We notice that the set $\F(\M, \bve, \s)$ appeared as
$\mathcal{N}^r(T, \mathbf{m}, A)$ in  Kasraoui \cite{Kasra08}, where
$\mathbf{m}$ is the column sum vector, and $A$ is the set of empty
rows, i.e., $A=\{i: \ve_i=0\}$.  Kasraoui also considered the set
$\mathcal{N}^c(T, \mathbf{n}, B)$ of  fillings whose row sum is an
arbitrary $\mathbb{N}$-vector $\mathbf{n}$ under the condition that
there is at most one $1$ in each column and where $B$ is the set of
empty columns. By a rotation of moon polyominoes, it is easily seen
that Theorem \ref{thm-row} and Theorem \ref{thm-column} also hold
for the set $\mathcal{N}^c(T, \mathbf{n}, B)$, as well as for the
set of fillings such that there is at most one $1$ in each row and
in each column.


\section{Mixed Statistics in Special Shapes}

In this section we show how Theorems \ref{thm-row} and
\ref{thm-column} specialize to classical combinatorial objects,
including  permutations, words, matchings, set partitions, and
simple graphs.

We first consider the  case of permutations and words. Fillings of
an $n \times m$ rectangle $\M$ are in bijection with words of length
$n$ on $[m]$. More precisely, a word $w=w_1w_2\cdots w_n$ on $[m]$
can be represented as a filling $M$ in which the cell in row $n+1-i$
and column $j$ is assigned the integer $1$
 if and only if $w_i=j$.
In the word $w_1w_2\cdots w_n$, a pair $(w_i, w_j)$ is an
\emph{inversion} if $i < j$ and $w_i > w_j$; we say that it is a
\emph{co-inversion}  if $i < j$ and $w_i < w_j$, see also
\cite{Mendes08}. Denote by $\inv(w)$ the  number of inversions of
$w$, and by $\asc(w)$ the number of co-inversions of $w$.

For $S \subseteq [n]$, the statistics $\alpha(S; M)$ and
$\beta(S;M)$ become
\begin{multline*}
 \hspace{1cm} \alpha(S; w) = \#\{ (w_i, w_j): n+1-j
 \in S \text{ and } (w_i, w_j) \text{ is a co-inversion}\} \\[3pt]
        + \#\{ (w_i, w_j) : n+1-j
         \not\in S \text{ and } (w_i, w_j) \text{ is an inversion}\}, \hspace{1cm}
\end{multline*}
and
\begin{multline*}
\hspace{1cm} \beta(S; w)= \#\{ (w_i, w_j): n+1-i \in S \text{ and }
(w_i, w_j) \text{ is a co-inversion}\} \\[3pt]
 + \#\{ (w_i, w_j) : n+1-i \not\in S \text{ and } (w_i, w_j) \text{ is an inversion}\}.\hspace{1cm}
\end{multline*}

For $T \subseteq [m]$, the statistics $\gamma(T; M)$  and
$\delta(T;M)$ become
\begin{multline*}
 \hspace{1cm}\gamma(T, w)= \#
  \{ (w_i, w_j): w_i \in T  \text{ and } (w_i, w_j) \text{ is a co-inversion}\}
  \\[3pt]
 + \#\{ (w_i, w_j) : w_j \not\in T \text{ and } (w_i, w_j) \text{ is an inversion}\},\hspace{1cm}
\end{multline*}
 and
\begin{multline*}
 \hspace{1cm}\delta(T, w)=
 \# \{ (w_i, w_j): w_j \in T  \text{ and } (w_i, w_j) \text{ is a co-inversion} \}
 \\[3pt]
 + \#\{ (w_i, w_j) : w_i \not\in T \text{ and } (w_i, w_j) \text{ is an inversion}\}.\hspace{1cm}
\end{multline*}

Let $W=\{1^{s_1}, 2^{s_2},\dots, m^{s_m}\}$ be a multiset with
$s_1+\cdots +s_m=n$. We adopt the  notation $R(W)$ for the set of
permutations, also called rearrangements,  of the elements in $W$.
Let $\lambda(A; w)$ denote any of the four statistics $ \alpha(S;
w), \beta(S; w), \gamma(T; w), \delta(T; w)$. Theorems \ref{thm-row}
and \ref{thm-column} imply that the bi-variate generating function
for $(\lambda(A; w), \lambda(\bar{A}; w))$ is symmetric and
\begin{eqnarray} \label{word}
\sum_{w \in R(W) } p^{\lambda(A; w)} q^{\lambda(\bar{A};w)} =
\sum_{w\in R(W)} p^{\inv(w)} q^{\asc(w)} =\genfrac[]{0pt}{}{n}{s_1,
\dots, s_m}_{p,q},
\end{eqnarray}
where $\genfrac[]{0pt}{}{n}{s_1, \dots, s_m}_{p,q}$ is the
$p,q$-Gaussian coefficient
\[
 \genfrac[]{0pt}{}{n}{s_1, \dots, s_m}_{p,q}
 = \frac{[n]_{p,q} ! }{ [s_1]_{p,q}! \cdots [s_m]_{p,q}! }.
\]
As usual, the $p,q$-integer $[r]_{p,q}$ is given by
\[
 [r]_{p,q} = \frac{p^r-q^r}{p-q} = p^{r-1} + p^{r-2}q + \cdots +pq^{r-2} +q^{r-1},
\]
and the $p,q$-factorial $[r]_{p,q}!$ is defined as $[r]_{p,q}!=
\prod_{i=1}^r [i]_{p,q}$.

We note that the symmetry of the distribution of $(\lambda(A; w),
\lambda(\bar{A}; w))$ can be easily seen from the map $w_1\cdots w_n
\rightarrow (m+1-w_1) \cdots (m+1-w_n)$ for $\alpha$ and $\beta$,
and the map $w_1\cdots w_n \rightarrow w_n \cdots w_1$ for $\gamma$
and $\delta$. Nevertheless, the generating function Eq.~\eqref{word}
seems to be new.  Chebikin \cite{Cheb} has considered the special
case of $\alpha(S; w)$ when $S$ is the set of even integers and $w$
ranges over all permutations of $[n]$.

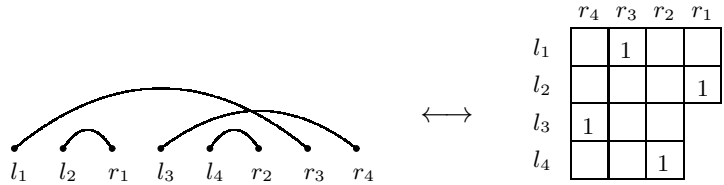
\begin{figure}[ht]
\begin{center}
\setlength{\unitlength}{0.5mm}
\begin{picture}(190,45)

\multiput(2,8)(13,0){8}{\circle*{1.5}}

\put(1,0){\footnotesize $l_1$} \put(14,0){\footnotesize $l_2$}

\put(27,0){\footnotesize $r_1$} \put(40,0){\footnotesize $l_3$}

\put(53,0){\footnotesize $l_4$} \put(65,0){\footnotesize $r_2$}

\put(79,0){\footnotesize $r_3$} \put(92,0){\footnotesize $r_4$}

\qbezier[500](2,8)(41,40)(80,8)   \qbezier[500](15,8)(21.5,18)(28,8)
\qbezier[500](41,8)(67,28)(93,8)  \qbezier[500](54,8)(60.5,18)(67,8)
\put(110,15){$\longleftrightarrow$}

\put(140,3){\footnotesize $l_4$} \put(140,13){\footnotesize $l_3$}

\put(140,23){\footnotesize $l_2$} \put(140,33){\footnotesize $l_1$}

\qbezier[500](150,0)(165,0)(180,0)
\qbezier[500](150,10)(165,10)(180,10)
\qbezier[500](150,20)(170,20)(190,20)
\qbezier[500](150,30)(170,30)(190,30)
\qbezier[500](150,40)(170,40)(190,40)
\qbezier[500](150,0)(150,20)(150,40)
\qbezier[500](160,0)(160,20)(160,40)
\qbezier[500](170,0)(170,20)(170,40)
\qbezier[500](180,0)(180,20)(180,40)
\qbezier[500](190,20)(190,30)(190,40)

\put(152,43){\footnotesize $r_4$} \put(162,43){\footnotesize $r_3$}

\put(172,43){\footnotesize $r_2$} \put(182,43){\footnotesize $r_1$}

\put(153,12){\footnotesize{1}} \put(163,32){\footnotesize{1}}
\put(173,2){\footnotesize{1}} \put(183.5,22){\footnotesize{1}}

\end{picture}
\end{center}

\caption{A matching and the corresponding filling of Ferrers
diagram.}\label{ferrers}

\end{figure}

We now consider the case of matchings and set partitions.
 As can be
seen in de Mier \cite{deMier} and Chen et al. \cite{Chen09},
general fillings of Ferrers diagrams correspond to multigraphs,
which include matchings, set partitions, and linked set
partitions. For simplicity, we give a description only for
matchings. Given a matching $\pi$ on $[2n]$, let $l_1 < l_2 <
\cdots < l_n$ be the left-hand endpoints, and $r_1 < r_2 < \cdots
< r_n$ be the right-hand endpoints. It determines  a Ferrers
diagram $\mathcal{F}$ whose rows are indexed by $l_1, \dots, l_n$
and columns are indexed by $r_n, \dots, r_1$, where a cell $(l_r,
r_k)$ is in the Ferrers diagram if and only if $l_r$ is on the
left of $r_k$. The cell $(l_r, r_k)$ is assigned the integer $1$
if and only if $(l_r, r_k)$ is an arc of the matching $\pi$. See
Figure \ref{ferrers} for an example.

A subset of rows corresponds to a subset $S$ of the left-hand
endpoints $\{ l_1, \dots, l_n\}$. The statistic $\alpha(S;M)$
corresponds to the mixed crossing-nesting statistic with respect to
the first left-hand endpoint. More precisely,  for a crossing formed
by two edges $(i_1,j_1)$ and $(i_2,j_2)$ with $i_1<i_2<j_1<j_2$, it
is said to be  an \emph{$S$-crossing} if $i_1 \in S$. Similarly, a
nesting formed by two edges $(i_1,j_1)$ and $(i_2,j_2)$ with
$i_1<i_2<j_2<j_1$ is said to be an \emph{ $S$-nesting } if $i_1 \in
S$. Thus the statistic $\alpha(S; M)$ becomes
\[
 \alpha(S; \pi) = \#\{\text{$S$-crossing of $\pi$}\} + \# \{\text{$\bar S$-nesting of $\pi$}\}.
\]
Theorem \ref{thm-row} asserts that $\sum_{\pi} p^{\alpha(S; \pi)}
q^{\alpha(\bar S; \pi)}$ is symmetric and independent of $S$, where
$\pi$ ranges over $P_n(A, B)$, the set of matchings with a given set
of left-hand endpoints $A$ and a given set of right-hand endpoints
$B$. In particular,  for each $r_i$, let
\[
h_i = \#\{\text{cell in the column indexed by $r_i$}\}-(i-1).
\]
By the generating function for the numbers
of crossings and nestings \cite{Chen08a,Kasra06}, we have
\begin{eqnarray} \label{cr-ne}
\sum_{\pi \in P_n(A, B)} p^{\alpha(S; \pi)} q^{\alpha(\bar S; \pi)}
 = \prod_{i=1}^n [h_i]_{p,q},
\end{eqnarray}
for any $S\subseteq \{l_1, \dots, l_n\}$. It is worth noting some
immediate consequences of Eq.~\eqref{cr-ne}. For example, for any
non-empty set $P_n(A, B)$, there is exactly one matching $\pi$ such
that $\alpha(S; \pi)=0$. It is not hard to construct such a
matching. Hence the number of matchings on $[2n]$ with $\alpha(S;
\pi)=0$ is given by the $n$-th Catalan number. Similar statements
hold when one considers the mixed crossing-nesting statistics with
respect to the second left-hand endpoint, the first right-hand
endpoint, and the second right-hand endpoint, respectively.

All the above results can be extended to set partitions
\cite{Kasra06} and linked set partitions \cite{Chen08b}, or more
generally, to simple graphs for which the left-degree of every
vertex is at most $1$, or the right-degree of every vertex is at
most $1$, see de Mier \cite{deMier}. Another way to see this is to
associate a simple graph with a filling of the triangular Ferrers
diagram $\Delta_n=(n-1, n-2, \dots, 1)$, see, for example,
\cite{Krattenthaler,deMier}.


\section{Proof of the Main Results}

It is sufficient to prove our results for $\alpha(S;M)$ and
$\gamma(T;M)$ only, since conclusions for $\beta(S; M)$ and
$\delta(T; M)$ can be obtained by reflecting  the moon polyomino
with respect to a horizontal line or a vertical line.

In Subsection 4.1, we recall Kasraoui's bijection $\Psi$ from
$\F(\M, \bve, \s)$ to   sequences of compositions \cite{Kasra08}.
Kasraoui's construction is stated for the set
$\mathcal{N}^c(T,\mathbf{n}, B)$. We shall give a description  to
fit our notation. The detailed justification of the bijection $\Psi$
can be found in \cite{Kasra08}, and hence is omitted. This bijection
will be used  in the proof of Lemma~\ref{lem-row} which states that
the pair  of the top-mixed statistics  $(\alpha(\{1\}; M),
\alpha(\overline{\{1\}}; M))$ is equidistributed with $(\se(M),
\ne(M))$. Theorem~\ref{thm-row} follows from an iteration of
Lemma~\ref{lem-row}. In Subsection 4.3 we provide two  proofs of
Theorem \ref{thm-column}. Again the crucial step is the observation
that $(\gamma(\{1\}; M), \gamma(\overline{\{1\}};M))$ has the same
distribution as $(\se(M), \ne(M))$.

\subsection{Kasraoui's bijection $\Psi$}

 If the columns of $\mathcal{M}$ are
$C_1,\ldots,C_m$ from left to right, it is clear that the sequence
of their lengths is unimodal and there exists a unique $k$ such that
\[
|C_1|\leq\cdots\leq|C_{k-1}|<|C_k|\geq|C_{k+1}|\geq\cdots\geq|C_m|,
\]
where $|C_i|$ is the length of the column $C_i$. The left part of
$\mathcal{M}$, denoted $L(\mathcal{M})$, is the set of columns
$C_i$'s with $1\leq i\leq k-1$, and the right part of
$\mathcal{M}$, denoted $R(\mathcal{M})$, is the set of columns
$C_i$'s with $k\leq i\leq m$. Note that the columns of maximal
length in $\mathcal{M}$ belong to $R(\mathcal{M})$.

We order the columns  $C_1,\ldots,C_m$ by a total order $\prec$ as
follows: $C_i\prec C_j$ if and only if
\begin{itemize}
  \item $|C_i|<|C_j|$ or
  \item $|C_i|=|C_j|$, $C_i\in L(\mathcal{M})$ and
        $C_j\in R(\mathcal{M})$, or
  \item $|C_i|=|C_j|$, $C_i,C_j\in L(\mathcal{M})$ and $C_i$ is on
        the left of $C_j$, or
  \item $|C_i|=|C_j|$, $C_i,C_j\in R(\mathcal{M})$ and $C_i$ is on
        the right of $C_j$.
\end{itemize}

For every column $C_i\in L(\mathcal{M})$, we define the rectangle
$\mathcal{M}(C_i)$ to be the largest rectangle that contains $C_i$
as the leftmost column. For $C_i\in R(\mathcal{M})$, the rectangle
$\mathcal{M}(C_i)$ is taken to be the largest rectangle that
contains $C_i$ as the rightmost column and does not contain any
column $C_j \in L(\mathcal{M})$ such that $C_j \prec C_i$.

Given $M\in\F(\mathcal {M},\bve,\s)$, we define a  coloring of $M$
by the following steps.

\noindent \textbf{The coloring of the filling $M$}
\begin{enumerate}
\item Color the cells of empty rows;
\item For each $C_i\in L(\mathcal{M})$, color the cells which
        are contained in the rectangle $\mathcal{M}(C_i)$ and on the
        right of any $1$-cell in $C_i$.
\item For each $C_i\in R(\mathcal{M})$, color the cells which
        are contained in the rectangle $\mathcal{M}(C_i)$ and on the
        left of any $1$-cell in $C_i$.
\end{enumerate}

Given $M$ with the coloring, let $\ce$ be a cell of $\M$. If $\ce$ is a $1$-cell,  we denote by $\auc(\ce;M)$ (resp. $\buc(\ce;M)$)
the number of uncolored empty cells in the same column as $\ce$ and above
(resp. below) $\ce$. If $\ce$ is empty,
we set $\auc(\ce;M)=\buc(\ce;M)=0$.
\begin{prop}\label{abuc}
Let $M\in\F(\mathcal {M},\bve, \s)$ and $\ce$ be a $1$-cell of
$C_i$.
\begin{enumerate}
  \item If $C_i\in L(\mathcal{M})$,  then $\auc(\ce;M)$
        \emph{(}resp. $\buc(\ce;M)$\emph{)} is equal to the
        number of NE \emph{(}resp. SE\emph{)} chains contained
        in the rectangle $\mathcal{M}(C_i)$ whose
        southwest \emph{(}resp. northwest\emph{)} $1$-cell is  $\ce$;
  \item If $C_i\in R(\mathcal{M})$, then $\auc(\ce;M)$
        \emph{(}resp. $\buc(\ce;M)$\emph{)} is equal to the
        number of SE \emph{(}resp. NE\emph{)} chains contained
        in the rectangle $\mathcal{M}(C_i)$ whose
        southeast \emph{(}resp. northeast\emph{)} $1$-cell is
        $\ce$.
\end{enumerate}
\end{prop}
\begin{exa}
\emph{Let $M$ be the 01-filling in Figure \ref{polyomino}, where
$L(\M)=\{C_1\}$ and $R(\M)=\{C_2, \dots, C_6\}$. Let $\ce$  be the
$1$-cell in the first column, and  $\ce'$ the $1$-cell in the fifth
column. Then $\auc(\ce; M)=1$, $\buc(\ce; M)=1$, $\auc(\ce'; M)=0$,
and $\buc(\ce'; M)=2$; See Figure \ref{auc-buc}.  }
\end{exa}

\begin{figure}[ht]
\begin{center}
\setlength{\unitlength}{0.4mm}
\begin{picture}(60,70)

\qbezier[500](10,0)(20,0)(30,0)   \qbezier[500](0,10)(30,10)(60,10)
\qbezier[500](0,20)(30,20)(60,20) \qbezier[500](0,30)(30,30)(60,30)
\qbezier[500](0,40)(30,40)(60,40) \qbezier[500](0,50)(25,50)(50,50)
\qbezier[500](0,60)(20,60)(40,60) \qbezier[500](10,70)(20,70)(30,70)
\qbezier[500](0,10)(0,35)(0,60)   \qbezier[500](10,0)(10,35)(10,70)
\qbezier[500](20,0)(20,35)(20,70) \qbezier[500](30,0)(30,35)(30,70)
\qbezier[500](40,10)(40,35)(40,60)\qbezier[500](50,10)(50,30)(50,50)
\qbezier[500](60,10)(60,25)(60,40)

\put(3,22.8){\footnotesize{1}} \put(13,2){\footnotesize{1}}
\put(23,12){\footnotesize{1}} \put(23,62){\footnotesize{1}}
\put(33,52){\footnotesize{1}} \put(43,32.8){\footnotesize{1}}

\put(5,25){\circle{8}} \put(45,35){\circle{8}}

\put(0,40){\line(1,1){10}} \put(0,45){\line(1,1){5}}
\put(5,40){\line(1,1){5}}

\put(0,30){\line(1,1){10}} \put(0,35){\line(1,1){5}}
\put(5,30){\line(1,1){5}}

\put(10,60){\line(1,1){10}} \put(10,65){\line(1,1){5}}
\put(15,60){\line(1,1){5}}

\put(10,50){\line(1,1){10}} \put(10,55){\line(1,1){5}}
\put(15,50){\line(1,1){5}}

\put(10,40){\line(1,1){10}} \put(10,45){\line(1,1){5}}
\put(15,40){\line(1,1){5}}

\put(10,30){\line(1,1){10}} \put(10,35){\line(1,1){5}}
\put(15,30){\line(1,1){5}}

\put(10,20){\line(1,1){10}} \put(10,25){\line(1,1){5}}
\put(15,20){\line(1,1){5}}

\put(10,10){\line(1,1){10}} \put(10,15){\line(1,1){5}}
\put(15,10){\line(1,1){5}}

\put(20,50){\line(1,1){10}} \put(20,55){\line(1,1){5}}
\put(25,50){\line(1,1){5}}

\put(20,40){\line(1,1){10}} \put(20,45){\line(1,1){5}}
\put(25,40){\line(1,1){5}}

\put(20,30){\line(1,1){10}} \put(20,35){\line(1,1){5}}
\put(25,30){\line(1,1){5}}

\put(20,20){\line(1,1){10}} \put(20,25){\line(1,1){5}}
\put(25,20){\line(1,1){5}}

\put(30,40){\line(1,1){10}} \put(30,45){\line(1,1){5}}
\put(35,40){\line(1,1){5}}

\put(30,30){\line(1,1){10}} \put(30,35){\line(1,1){5}}
\put(35,30){\line(1,1){5}}

\put(30,20){\line(1,1){10}} \put(30,25){\line(1,1){5}}
\put(35,20){\line(1,1){5}}

\put(40,40){\line(1,1){10}} \put(40,45){\line(1,1){5}}
\put(45,40){\line(1,1){5}}

\end{picture}
\end{center}

\caption{The statistics $\mathrm{auc}$ and $\mathrm{buc}$ for cells
in a 01-filling.} \label{auc-buc}

\end{figure}

The following theorem can be deduced from  Prop.~\ref{abuc}.

\begin{thm}\label{thmd}
\[\ne(M)=\sum\limits_{\ce\in L(\mathcal{M})}\auc(\ce;M)+
\sum\limits_{\ce\in R(\mathcal{M})}\buc(\ce;M),\]
\[\se(M)=\sum\limits_{\ce\in L(\mathcal{M})}\buc(\ce;M)+
\sum\limits_{\ce\in R(\mathcal{M})}\auc(\ce;M).\]
\end{thm}

For $M \in \F(\mathcal {M},\bve, \s)$, let $a_i$ be the number of
empty rows (i.e., $\{R_i: \ve_i=0\}$) that intersect the column
$C_i$. Suppose that $C_{i_1}\prec C_{i_2}\prec\cdots\prec C_{i_m}$.
For $j=1,\ldots,m$, we define
\begin{eqnarray} \label{h_i}
h_{i_j}=|C_{i_j}|-a_{i_j}-(s_{i_1}+s_{i_2}+\cdots+s_{i_{j-1}}).
\end{eqnarray}
Note that the numbers $h_i$'s have  the following interpretation. If
one puts $1$-cells in the columns of $M$ from the smallest to the
largest under the order $\prec$, then $h_{i_j}$ is the number of
available cells in the $j$-th column to be filled. For positive
integers $n$ and $k$, denote by $\mathcal{C}_k(n)$ the set of
compositions of $n$ into $k$ nonnegative parts, that is,
$\mathcal{C}_k(n)=\{(\lambda_1,\lambda_2,\ldots,\lambda_k)\in
\mathbb{N}^k: \sum_{i=1}^k \lambda_i=n \}$. The bijection $\Psi$ is
constructed as follows.

\noindent \textbf{The bijection} $\Psi: \F(\mathcal {M},\bve,
\s)\longrightarrow \mathcal{C}_{s_1+1}(h_1-s_1)\times
\mathcal{C}_{s_2+1}(h_2-s_2) \times \cdots \times
\mathcal{C}_{s_m+1}(h_m-s_m).$

 For each $M\in\F(\mathcal
{M},\bve,\s)$ with the coloring, $\Psi(M)$ is a sequence of
compositions $(c^{(1)},c^{(2)},\ldots,c^{(m)})$, where
\begin{itemize}
 \item $c^{(i)}=(0)$ if $s_i=0$. Otherwise
 \item $c^{(i)}= (c^{(i)}_{1},c^{(i)}_{2},\ldots,c^{(i)}_{s_i+1})$ where
  \begin{itemize}
  \item $c^{(i)}_1$ is the number of uncolored cells above the first
        $1$-cell in the column $C_i$;
  \item $c^{(i)}_k$ is the number of uncolored cells between the
        $(k-1)$-st and the $k$-th $1$-cells in the column $C_i$, for $2 \leq k \leq s_i$;
  \item $c^{(i)}_{s_i+1}$ is the number of uncolored cells below
        the last $1$-cell in the column $C_i$.
\end{itemize}
\end{itemize}

Let $\mathbf{c}=\Psi(M)=(c^{(1)},c^{(2)},\ldots,c^{(m)})$, and
 $\ce$ be the $k$-th $1$-cell in the column $C_i$. It
follows from the bijection $\Psi$ that
\[\auc(\ce;M)=c^{(i)}_1+c^{(i)}_2+\cdots+c^{(i)}_k,\]
\[\buc(\ce;M)=c^{(i)}_{k+1}+c^{(i)}_{k+2}+\cdots+c^{(i)}_{s_i+1}
=h_i-s_i-(c^{(i)}_1+c^{(i)}_2+\cdots+c^{(i)}_k).\]

Now Theorem \ref{thmd} can be rewritten as

\begin{thm}\label{thma}
Let $M \in \F(\mathcal {M},\bve,\s)$ and
$\mathbf{c}=\Psi(M)= (c^{(1)},c^{(2)},\ldots,c^{(m)})$. Then
\[\ne(M)=\sum\limits_{C_i\in L(\mathcal{M})}
\sum\limits_{k=1}^{s_i}(c_1^{(i)}+c_2^{(i)}+\cdots+c_k^{(i)}) +
\sum\limits_{C_j\in R(\mathcal{M})}\sum\limits_{k=1}^{s_j}
(h_j-s_j-c_1^{(j)}-c_2^{(j)}-\cdots-c_k^{(j)}),\]
\[\se(M)=\sum\limits_{C_i\in L(\mathcal{M})}
\sum\limits_{k=1}^{s_i}(h_i-s_i-c_1^{(i)}-c_2^{(i)}-\cdots-c_k^{(i)})
+ \sum\limits_{C_j\in R(\mathcal{M})}\sum\limits_{k=1}^{s_j}
(c_1^{(j)}+c_2^{(j)}+\cdots+c_k^{(j)}).\]
\end{thm}

Summing over the sequences of compositions yields the symmetric
generating function.

\begin{thm}[Kasraoui] \label{gf}
 \[\sum_{M \in \F(\M, \bve,\s)} p^{\ne(M)} q^{\se(M)} =\sum_{M \in \F(\M, \bve,\s)} p^{\se(M)} q^{\ne(M)}
=\prod_{i=1}^m  \genfrac[]{0pt}{}{h_i}{s_i}_{p,q}.\]
\end{thm}

\subsection{Proof of Theorem \ref{thm-row}}

To prove Theorem \ref{thm-row} for the top-mixed  statistic
$\alpha(S;M)$, we first consider the  special case when
$\mathcal{R}(S)$ contains the first row only.

\begin{lem}\label{lem-row}
For $S=\{1\}$, we have
\[
F^t_{\{1\}}(p,q)=F^t_{\emptyset}(p,q)=
\sum\limits_{M\in\F(\mathcal{M},\bve,\s)}p^{\se(M)}q^{\ne(M)}.\]
\end{lem}
\pf We assume that the first row is nonempty. Otherwise the identity
is obvious. Given a filling $M\in\F(\mathcal{M},\bve,\s)$, assume
that the unique $1$-cell  of the first row lies in the column $C_t$.
Let the upper polyomino $\M_u$ be the union of the rows that
intersect $C_t$, and the lower polyomino $\M_d$ be the complement of
$\M_u$, i.e., $\M_d=\M \setminus \M_u$. We aim to construct a
bijection $\phi_\alpha: \F(\M,\bve, \s) \rightarrow \F(\M,\bve, \s)$
such that for any filling $M$,
\[
(\alpha(\{1\};M), \alpha(\overline{\{1\}}; M) ) = (\se(\phi_\alpha(M)), \ne(\phi_\alpha(M))),
\]
and $\phi_\alpha(M)$ is identical to $M$ on  $\M_d$.

Let $M_u = M \cap \M_u$ and $M_d=M \cap \M_d$. Let $s_i'$ be the
number of $1$-cells of $M$ in the column $C_i \cap \M_u$, and
$\s'=(s_1', \dots, s_m')$. Let $\bve'=(\ve_1, \dots, \ve_r)$, where
$r$ is the number of rows in $\M_u$. We shall define  $\phi_\alpha$
on $\F(\M_u, \bve',\s')$  such that $\phi_\alpha(M_u) \in \F(\M_u,
\bve', \s')$  and
\[
(\alpha(\{1\};M_u), \alpha(\overline{\{1\}}; M_u) ) = (\se(\phi_\alpha(M_u)), \ne(\phi_\alpha(M_u))).
\]

Let $C_i'=C_i \cap \M_u$. Suppose that in $\M$ the columns
intersecting with the first row are  $C_a,\ldots,C_t,\ldots,C_b$
from left to right. Then $C_t=C_t'$, and in $\M_u$ the columns
$C_a', \ldots,C'_t,\ldots,C'_b$  intersect the first row. Assume
that among them the ones with the same length as $C'_t$ are
$C'_u,\ldots,C'_t,\ldots,C'_v$ from left to right. Clearly, the
columns $C'_u,\ldots,C'_t,\ldots,C'_v$ are those with maximal length
and belong to $R(\M_u)$. Note that in $M_u$, the number of top
$\{1\}$-NE chains is $\sum_{a\leq i<t}s_i'$,  while the number of
top  $\{1\}$-SE chains is  $\sum_{t<i\leq b}s_i'$. Let $h_i'$ be
given as in Eq.~\eqref{h_i} for $\F(\M_u, \bve', \s')$. Let
$\textbf{c}=\Psi(M_u)= (c^{(1)},c^{(2)},\ldots,c^{(m)})$, from
Theorem \ref{thma} we see that \allowdisplaybreaks
\begin{eqnarray}
\alpha(\{1\};M_u)
&=& \sum\limits_{a\leq i<t}s_i'+\sum\limits_{C'_i\in
L(\M_u)} \sum\limits_{k=1}^{s_i'}(h_i'-s_i'-c_1^{(i)}-c_2^{(i)}
-\cdots-c_k^{(i)}) \nonumber \\
&&{}+\sum\limits_{C'_j\in R(\M_u)}\sum\limits_{k=1}^{s_j'}
(c_1^{(j)}+c_2^{(j)}+\cdots+c_k^{(j)})-\sum\limits_{t<i\leq b}s_i' \nonumber \\
&=& \sum\limits_{a\leq i<u}s_i'+(h_t'-s_t')+\sum\limits_{C'_i\in L(\M_u)}
\sum\limits_{k=1}^{s_i'}(h_i'-s_i'-c_1^{(i)}-c_2^{(i)}
-\cdots-c_k^{(i)}) \nonumber \\
&&{}+\sum\limits_{C'_j\in R(\M_u)}\sum\limits_{k=1}^{s_j'}
(c_1^{(j)}+c_2^{(j)}+\cdots+c_k^{(j)})-\sum\limits_{t<i\leq b}s_i'. \label{alpha_1}
\end{eqnarray}
The second equation holds since $C'_t \prec C'_{t-1} \prec  \cdots
\prec C'_u$ are the largest $t-u+1$ columns in $R(\M_u)$ under the
order $\prec$. By definition $h_t'$ is the number of available rows
when all the smaller columns of $\M_u$ have been filled. Those
available rows will be filled by the $1$'s in the columns $C'_t,
\dots, C'_u$. Hence $h_t'=s_t'+\cdots +s_u'$. Similarly, we have
\begin{eqnarray}
\alpha(\overline{\{1\}};M_u)&=&\sum\limits_{t<i\leq b}s_i'+\sum\limits_{C'_i\in
L(\M_u)}
\sum\limits_{k=1}^{s_i'}(c_1^{(i)}+c_2^{(i)}+\cdots+c_k^{(i)}) \nonumber\\
&&{}+\sum\limits_{C'_j\in R(\M_u)}\sum\limits_{k=1}^{s_j'}
(h_j'-s_j'-c_1^{(j)}-c_2^{(j)}-\cdots-c_k^{(j)})-\sum\limits_{a\leq
i<t}s_i'\nonumber \\
&=&\sum\limits_{t<i\leq b}s_i'+\sum\limits_{C'_i\in L(\M_u)}
\sum\limits_{k=1}^{s_i'}(c_1^{(i)}+c_2^{(i)}+\cdots+c_k^{(i)}) \nonumber \\
&&{}+\sum\limits_{C'_j\in R(\M_u)}\sum\limits_{k=1}^{s_j'}
(h_j'-s_j'-c_1^{(j)}-c_2^{(j)}-\cdots-c_k^{(j)})-\sum\limits_{a\leq
i<u}s_i'-(h_t'-s_t'). \label{alpha_2}
\end{eqnarray}

The fact that the $1$-cell of the first row lies in the column
$C'_t$ implies that
 $c^{(t)}_1=0$, and $c^{(i)}_1 > 0$ for $a\leq i<u$ or $t<i\leq b$.
We define the filling $\phi_\alpha(M_u)$ by setting
$\phi_\alpha(M_u) = \Psi^{-1}(\tilde{\textbf{c}})$, where
$\tilde{\textbf{c}}$  is obtained from $\mathbf{c}$ as follows:
\begin{eqnarray*}
\left\{ \begin{array}{ll}
 \tilde{c}^{(i)}=
        (c^{(i)}_1-1,c^{(i)}_2,\ldots,c^{(i)}_{s_i},
        c^{(i)}_{s_i+1}+1), & \text{if   $a\leq i<u$ or $t<i\leq b$, and $s_i'\neq
        0$},\\[5pt]
  \tilde{c}^{(t)}=(c^{(t)}_2,c^{(t)}_3,
        \ldots,c^{(t)}_{s_t+1},c^{(t)}_1), & \text{if $i=t$},\\[5pt]
   \tilde{c}^{(i)}=c^{(i)}, & \text{for any other $i$}.
\end{array} \right.
\end{eqnarray*}

Comparing the formulas \eqref{alpha_1} and \eqref{alpha_2} with
Theorem \ref{thma} for $\tilde{\textbf{c}}$, it is  easily verified
that
\[
(\alpha(\{1\};M_u),\alpha(\overline{\{1\}};M_u))=(\se(\phi_\alpha(M_u)),\ne(\phi_\alpha(M_u))).
\]
 Now  $\phi_\alpha(M)$ is obtained from $M$ by replacing $M_u$ with  $\phi_\alpha(M_u)$.

\noindent {Claim}: $(\alpha(\{1\};M), \alpha(\overline{\{1\}}; M) )
= (\se(\phi_\alpha(M)), \ne(\phi_\alpha(M)))$ for any $M \in \F(\M,
\bve, \s)$.

This is true because (1) $M$ has the same number of top
$\{1\}$-NE/SE chains as $M_u$, since every top $\{1\}$-NE/SE chains
of $M$ must appear in $M_u$; (2) $M_d$ appears in both $M$ and
$\phi_\alpha(M)$; (3) If $(\ce, \ce')$ is an NE chain or an SE chain
with $\ce \in \M_u$ and $\ce' \in \M_d$, by the intersection-free
property of $\M$, both $\ce$ and $\ce'$ are in columns $\{C_a,
\dots, C_b\}$. For any fixed $\ce' \in \M_d$, the number of NE
(resp. SE) chains formed by $\ce'$ and $1$-cells $\ce$ in the column
$C_j \cap \M_u$ is unchanged under the map $\phi_\alpha$ since
$\phi_\alpha$ preserves the column sum and row sum of $M_u$.

To show that  $\phi$ is a bijection on $\mathbf{F}(\mathcal
{M},\bve, \s)$, it is enough to explain how to determine from
$\phi_\alpha(M)$ the column $C_t$, and hence the upper polyomino
$\M_u$.  Then the correspondence between $\mathbf{c}$ and
$\tilde{\mathbf{c}}$ becomes obvious. To this end,  we shall use the
map $\Psi$ defined in Subsection 4.1. If the columns intersecting
the first row are $C_a, \dots, C_b$ in $\M$, then $C_t$ is the
smallest column in $\{C_a, \dots, C_b\}$ under the order $\prec$
with the property that the last entry $c^{(t)}_{s_t+1}$ is $0$  in
the corresponding composition $c^{(t)}$. The rest of the proof is
straightforward. \qed

\begin{prop}\label{prop-row}
Assume $S=\{r_1,r_2,\ldots, r_s\} \subseteq [n]$ with
$r_1<r_2<\cdots<r_s$. Let $S'=\{r_1, r_2,\ldots, r_{s-1}\}$. Then
$F^t_S(p,q)=F^t_{S'}(p,q)$.
\end{prop}
\pf Let $X=\{R_i: 1 \leq i < r_s\}$ be the set of rows above the row $R_{r_s}$,
 and $Y$ be the set of remaining rows.
Given a filling $M\in \textbf{F}(\mathcal {M},\bve,\s)$, let
$\mathcal{T}(M)$ be the set of fillings $M'\in \textbf{F}(\mathcal
{M},\bve,\s)$ that are identical to $M$ in the rows of $X$.
Construct a bijection $\theta_{r_s}\colon  \mathcal{T}(M)
\rightarrow \mathcal{T}(M)$ by setting $\theta_{r_s}(M)$ to be the
filling obtained from $M$ by replacing $M \cap Y$ with
$\phi_\alpha(M \cap Y)$.

\begin{figure}[ht]
\begin{center}
\setlength{\unitlength}{0.35mm}
\begin{picture}(170,160)

\qbezier[500](70,0)(75,0)(80,0) \qbezier[500](60,10)(70,10)(80,10)
\qbezier[500](60,20)(70,20)(80,20)
\qbezier[500](50,30)(75,30)(100,30)
\qbezier[500](80,40)(95,40)(110,40)
\qbezier[500](30,50)(40,50)(50,50)
\qbezier[500](80,50)(100,50)(120,50)
\qbezier[500](20,60)(30,60)(40,60)
\qbezier[500](80,60)(100,60)(120,60)
\qbezier[500](10,70)(25,70)(40,70)
\qbezier[500](80,70)(105,70)(130,70)
\qbezier[500](0,80)(20,80)(40,80)
\qbezier[500](80,80)(105,80)(130,80)
\qbezier[500](0,90)(20,90)(40,90)
\qbezier[500](80,90)(110,90)(140,90)
\qbezier[500](0,100)(20,100)(40,100)
\qbezier[500](80,100)(110,100)(140,100)
\qbezier[500](0,110)(20,110)(40,110)
\qbezier[500](80,110)(110,110)(140,110)

\qbezier[500](20,120)(70,120)(120,120)
\qbezier[500](30,130)(70,130)(110,130)
\qbezier[500](40,140)(75,140)(110,140)
\qbezier[500](40,150)(75,150)(110,150)

\qbezier[500](0,80)(0,95)(0,110) \qbezier[500](10,70)(10,90)(10,110)
\qbezier[500](20,60)(20,90)(20,120)
\qbezier[500](30,50)(30,90)(30,130)
\qbezier[500](40,50)(40,100)(40,150)
\qbezier[500](50,30)(50,40)(50,50)
\qbezier[500](50,120)(50,135)(50,150)
\qbezier[500](60,10)(60,20)(60,30)
\qbezier[500](60,120)(60,135)(60,150)
\qbezier[500](70,0)(70,15)(70,30)
\qbezier[500](70,120)(70,135)(70,150)
\qbezier[500](80,0)(80,75)(80,150)
\qbezier[500](90,30)(90,90)(90,150)
\qbezier[500](100,30)(100,90)(100,150)
\qbezier[500](110,40)(110,95)(110,150)
\qbezier[500](120,50)(120,85)(120,120)
\qbezier[500](130,70)(130,90)(130,110)
\qbezier[500](140,90)(140,100)(140,110)

\put(81.5,133){\footnotesize ce} \put(59,73){\footnotesize ce'}
\put(45,152){\vector(0,1){10}}   \put(40,165){\small$C_a$}
\put(85,152){\vector(0,1){10}}   \put(80,165){\small$C_t$}
\put(105,152){\vector(0,1){10}} \put(100,165){\small$C_b$}
\put(125,115){\vector(1,0){20}} \put(147,111){\small$R_{r_s}$}
\put(106,35){\vector(1,0){20}} \put(127,32){\small$R_p$}
\put(118,135){\vector(1,0){20}} \put(140,130){\small$R_u$}

\linethickness{1.5pt}

\put(20,120){\line(1,0){100}}

\multiput(-30,119.8)(5,0){10}{\line(0,1){1}}

\thinlines

\put(-25,140){\Large $X$}    \put(-25,60){\Large $Y$}

\put(20,110){\line(1,1){10}} \put(20,115){\line(1,1){5}}
\put(25,110){\line(1,1){5}}

\put(30,110){\line(1,1){10}} \put(30,115){\line(1,1){5}}
\put(35,110){\line(1,1){5}}

\put(80,140){\line(1,1){10}} \put(80,145){\line(1,1){5}}
\put(85,140){\line(1,1){5}}

\put(80,120){\line(1,1){10}} \put(80,125){\line(1,1){5}}
\put(85,120){\line(1,1){5}}

\put(80,110){\line(1,1){10}} \put(80,115){\line(1,1){5}}
\put(85,110){\line(1,1){5}}

\put(80,100){\line(1,1){10}} \put(80,105){\line(1,1){5}}
\put(85,100){\line(1,1){5}}

\put(80,90){\line(1,1){10}} \put(80,95){\line(1,1){5}}
\put(85,90){\line(1,1){5}}

\put(80,80){\line(1,1){10}} \put(80,85){\line(1,1){5}}
\put(85,80){\line(1,1){5}}

\put(80,70){\line(1,1){10}} \put(80,75){\line(1,1){5}}
\put(85,70){\line(1,1){5}}

\put(80,60){\line(1,1){10}} \put(80,65){\line(1,1){5}}
\put(85,60){\line(1,1){5}}

\put(80,50){\line(1,1){10}} \put(80,55){\line(1,1){5}}
\put(85,50){\line(1,1){5}}

\put(80,40){\line(1,1){10}} \put(80,45){\line(1,1){5}}
\put(85,40){\line(1,1){5}}

\put(80,30){\line(1,1){10}} \put(80,35){\line(1,1){5}}
\put(85,30){\line(1,1){5}}

\put(90,110){\line(1,1){10}} \put(90,115){\line(1,1){5}}
\put(95,110){\line(1,1){5}}

\put(100,110){\line(1,1){10}} \put(100,115){\line(1,1){5}}
\put(105,110){\line(1,1){5}}

\put(110,110){\line(1,1){10}} \put(110,115){\line(1,1){5}}
\put(115,110){\line(1,1){5}}

\end{picture}
\end{center}

\caption{NE chains formed by the cell $\mathrm{ce} \in X$ and
$Y$-cells.} \label{proof}

\end{figure}
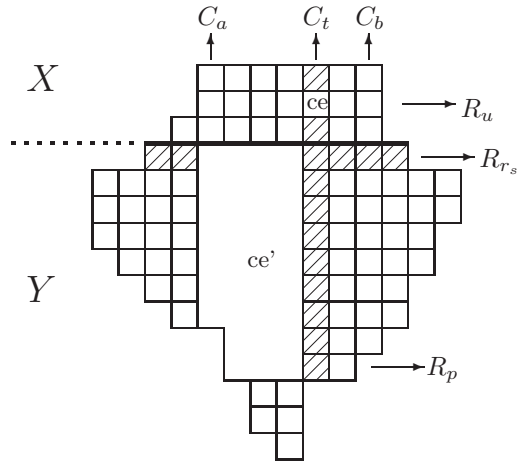

We proceed to show that
\begin{eqnarray} \label{eqn1}
(\alpha(S;M),\alpha(\bar S;M)) =
(\alpha(S'; \theta_{r_s}(M)), \alpha(\overline{S'}; \theta_{r_s}(M))).
\end{eqnarray}
There are three cases.
\begin{description}
\item [Case 1]  An NE or an SE chain consisting of two $1$-cells in $X$ contributes equally to both pairs of statistics.

\item [Case 2] By Lemma~\ref{lem-row}, the set of NE chains and SE chains consisting of two cells in $Y$
contributes equally to both pairs of statistics.

\item [Case 3] For a $1$-cell $\ce$ in $X$, assume $\ce$ is in row $R_u$ and  column $C_t$.
Let $T=\{C_a, \dots, C_b\}$ be the set of columns intersecting both
the rows $R_{r_s}$ and $R_u$, and  $R_p$ ($p \geq r_s$) be the
lowest row that intersects $C_t$. If $\ce$ forms an NE chain with a
cell $\ce'$ in $Y$, then $\ce'$ is in a row on or above $R_p$, and
in a column in $\{C_a, \dots, C_{t-1}\}$.

It follows that the number of NE chains of the form $(\ce, \ce')$
for a fixed $1$-cell $\ce \in X$ equals the number of $1$-cells in
the area $\{ (R_i, C_j): r_s \leq i \leq p, a \leq j < t\}$, see
Figure \ref{proof}.
  This number is unchanged under the map $\phi_\alpha$, as $\phi_\alpha$
preserves the column sum and the row sum, and hence the number of
$1$'s in columns $C_a, \dots, C_{t-1}$, and the number of $1$'s in
rows $\{R_i:i > p\}$. Similarly, the number of SE chains $(\ce,
\ce')$ with $\ce \in X$ and $\ce' \in Y$ is unchanged under the map
$\phi_\alpha$. Thus NE and SE chains formed by one $X$-cell and one
$Y$-cell contribute equally to the two pairs of statistics as well.
\end{description}

Thus \eqref{eqn1} is proved by combining the above three cases. \qed

\noindent \emph{Proof of Theorem \ref{thm-row}}. Assume $S=\{r_1,
r_2,\ldots, r_s\} \subseteq \mathcal{R}$ with $r_1<r_2<\cdots<r_s$.
Let $\Theta_\alpha  = \theta_{r_1} \circ \theta_{r_2} \circ \cdots
\circ \theta_{r_s}$, where $\theta_r$ is defined in the proof of
Prop.~\ref{prop-row}. Then $\Theta_\alpha$ is a bijection on $\F(\M,
\bve, \s)$
 with the property that
\[
 (\alpha(S; M), \alpha(\bar S; M)) = (\se(\Theta_\alpha(M)), \ne(\Theta_\alpha(M))).
\]
The symmetry of $F^t_S(p,q)$ follows from Theorem \ref{gf}.
\qed

\subsection{Proof of Theorem \ref{thm-column} }

 Theorem \ref{thm-column} is concerned with  the left-mixed statistic
$\gamma(T;M)$. The proof is similar to that of Theorem
\ref{thm-row}. The key idea amounts to the observation that
$(\gamma(\{1\};M),\gamma(\overline{\{1\}};M))$ is equidistributed
with $(\se(M), \ne(M))$. We provide two proofs of this fact: one
is based on generating functions, and the other is bijective.

\begin{lem}\label{lem-column}
For $T=\{1\}$, we have
\[
  G^l_{\{1\}}(p,q) = G^l_\emptyset(p,q)= \prod_{i=1}^m \genfrac[]{0pt}{}{h_i}{s_i}_{p,q}.
\]
\end{lem}
\noindent \emph{First proof of Lemma \ref{lem-column}.} We conduct
induction on the number of columns of $\mathcal{M}$. The statement
is trivial if $\mathcal{M}$ has only one column.

Assume that Lemma \ref{lem-column} holds for 01-fillings on any moon
polyominoes with less than $m$ columns. Suppose that  $\mathcal{M}$
have $m$ columns. Consider the minimal column $C$ under the order
$\prec$. There are two cases.
\begin{enumerate}
\item  $C=C_1$ is the leftmost column of $\mathcal{M}$.
In this case we employ the bijection $\Psi$. For any filling $M$
with $\Psi(M)= (c^{(1)}, c^{(2)},\dots, c^{(m)})$, let $\tau(M) =
\Psi^{-1}(c^{(1),r}, c^{(2)}, \ldots, c^{(m)})$, where
\[
c^{(1),r} = (c^{(1)}_{s_1+1}, \dots, c^{(1)}_2,c^{(1)}_1) \text{ if }
c^{(1)}=(c^{(1)}_1, c^{(1)}_2, \dots, c^{(1)}_{s_1+1}).
\]
It is readily checked that $(\gamma(\{1\}; M),
\gamma(\overline{\{1\}}; M) ) = (\se(\tau(M)), \ne(\tau(M)))$.

\item $C=C_m$ is the rightmost column of $\mathcal{M}$. We first
prove the case for rectangular shapes. Assume $\mathcal{M}$ is a
rectangle  with $n$ non-empty rows. A filling $M$ of $\mathcal{M}$
can be read as a word $w=w_1w_2\cdots w_ n$ where $w_i=j$ if the
only 1-cell in the  $(n+1-i)$-th non-empty row appears in the
$j$-th column. It is clear that $ \se(M)= \text{inv}(w_1w_2\cdots
w_n)$ and
 $\ne(M) =\text{coinv}(w_1w_2\cdots w_n)$.
In addition, fillings with a given column sum $\s=(s_1, \dots, s_m)$
correspond to words on the multiset $\{1^{s_1}, \dots, m^{s_m}\}$.
Therefore
\begin{eqnarray} \label{One}
 \sum_{M \in  \mathbf{F}(\mathcal {M},\bve, \mathbf{s})} p^{\se(M)}q^{\ne(M)} =
\genfrac[]{0pt}{}{n}{s_1, \dots, s_m}_{p,q}.
\end{eqnarray}

Observe that
\[
 \gamma(\{1\}; M) = \#~\{(w_i, w_j)\ | \ w_i=1 <  w_j \text{ and } i < j \}
 + \#~\{(w_i, w_j)\ | \ i < j \text{ and } w_i > w_j \neq 1 \}.
\]
Let $\epsilon(w_1w_2\cdots w_n)=\epsilon(w_1)\epsilon(w_2) \cdots \epsilon(w_n)$
where
\begin{eqnarray*}
 \epsilon(w_i) = \left\{ \begin{array}{ll}
 m+1, & \text{ if } w_i=1,\\[3pt]
 w_i, & \text{ otherwise}. \end{array} \right.
\end{eqnarray*}

Then $ \gamma(\{1\}; M) = \text{inv}(\epsilon(w_1w_2\cdots  w_n))$.
Similarly,  $\gamma(\overline{\{1\}}; M) =
\asc(\epsilon(w_1w_2\cdots w_n)) $. When $M$ ranges over $\F(\M,
\bve, \s)$, the word $\epsilon(w_1w_2\cdots w_n)$ ranges over all
rearrangements of the multiset $W= \{(m+1)^{s_1}, 2^{s_2}, \dots,
m^{s_m}\}$. Hence
\begin{eqnarray}\label{Two}
 \sum_{M \in  \mathbf{F}(\mathcal {M},\bve, \mathbf{s})} p^{\gamma(\{1\};M)}q^{\gamma(\overline{\{1\}};M
)} = \sum_{w\in R(W)} p^{\inv(w)}q^{\asc(w)}
 =\genfrac[]{0pt}{}{n}{s_1, \dots, s_m}_{p,q}.
\end{eqnarray}
Comparing \eqref{One} and \eqref{Two}, we complete the proof for a
rectangular shape $\mathcal{M}$. \vspace{.2cm}

Now we deal with the case for a general shape $\mathcal{M}$. Let
$\mathcal{M}(C) $ be the largest rectangle that contains $C$. Let
$\mathcal{M}_1 = \mathcal{M} \setminus C$ and $\mathcal{M}_1(C) =
\mathcal{M}(C)\setminus C $. By the inductive hypothesis,
$(\gamma(\{1\}; M), \gamma(\overline{\{1\}}; M))$ has the same
distribution as $(\se(M), \ne(M))$ over the set of fillings in
$\mathbf{F} (\mathcal {M}_1,\bve', \mathbf{s}-\{s_m\})$, where
$\bve'$ is  the row-vector when one removes the rows with a $1$-cell
in the column $C$. We analyze the  contribution when one adds the
last column $C$ with $s_m$ many 1-cells.

Given a filling $M$ on $\mathcal{M}_1$, let $\mathcal{S}(M)$ be the set of all fillings that
have
the same 1-cells in rows other than those in $\mathcal{M}_1(C)$. For any $N \in  \mathcal{S}(
M)$,
let $N(C)$ be the restriction of $N$ on $\mathcal{M}_1(C)$, then
the values
\[
\gamma(\{1\}; N)- \gamma(\{1\}; N(C)) \quad \text{and} \quad \gamma(\overline{\{1\}}; N)- \gamma(\overline{\{1\}}; N(C)) \
\]
are  constants over $\mathcal{S}(M)$, which will not change when the
last column $C$ is added. Note that $\mathcal{M}_1(C)$ is a
rectangular shape. Hence $N(C)$ can be identified as a word of
length $n$ on $\{1^{s_1'}, 2^{s_2'}, \dots, (m-1)^{s_{m-1}'}\}$,
where $n$ is the number of non-empty rows in $\mathcal {M}_1(C)$,
and
 $s_i'$ is the number of $1$'s of $M$ in $C_i \cap \mathcal{M}(C)$.
By the argument for a rectangular shape, we find that
\[
 \sum_{N \in \mathcal{S}(M)} p^{\gamma(\{1\}; N(C))} q^{\gamma(\overline{\{1\}}; N(C))} =
\genfrac[]{0pt}{}{n}{s_1', \dots, s_{m-1}'}_{p,q}.
\]

Adding the last column $C$ with $s_m$ $1$-cells is equivalent to
inserting $s_m$ many $m$'s to a word on $\{1^{s_1'}, 2^{s_2'}, \dots, (m-1)^{s_{m-1}'}\}$.
Again using the transformation $\epsilon(w_1w_2 \cdots w_{n+s_m})=
\epsilon(w_1)\epsilon(w_2)\cdots \epsilon(w_{n+s_m})$,
and assuming $W'=\{ (m+1)^{s_1'}, 2^{s_2'}, \dots, (m-1)^{s_{m-1}'}, m^{s_m}\}$,
we have
\begin{eqnarray*}
\sum_{M \in \mathbf{F}(\mathcal{M}(C),\bve_1, \mathbf{s'}
\cup\{s_m\})} p^{\gamma(\{1\}; M)} q^{\gamma(\overline{\{1\}}; M)}
&= &\sum_{w\in R(W') }
 p^{\inv(w)} q^{\asc(w)}  \\
&= &\genfrac[]{0pt}{}{n+s_m}{s_1', \dots, s_{m-1}', s_m}_{p,q},
\end{eqnarray*}
where $\bve_1$ is the restriction of $\bve$ on the rows in $\M(C)$.
Thus we deduce that the contribution of the last column $C$ over the
set $\mathcal{S}(M)$ is given by
\[\genfrac[]{0pt}{}{n+s_m}{s_1', \dots, s_{m-1}', s_m}_{p,q}
\Big/
\genfrac[]{0pt}{}{n}{s_1', \dots, s_{m-1}'}_{p,q}
=\genfrac[]{0pt}{}{n+s_m}{s_m}_{p,q},
\]which is independent of $s_1', \dots, s_{m-1}'$.

Summing over all distinct sets of the form $\mathcal{S}(M)$, we
conclude that adding the last column $C$ contributes a factor of $
\genfrac[]{0pt}{}{n+s_m}{s_m}_{p,q}$ to $G^l_{\{1\}}(p,q)$. It
follows from \eqref{h_i} that  $n+s_m = h_m$. Hence an inductive
argument yields
\begin{equation*}
 G^l_{\{1\}}(p,q)= \prod_{i=1}^m \genfrac[]{0pt}{}{h_i}{s_i}_{p,q}
= G^l_{\emptyset} (p,q).\tag*{\rule{4pt}{7pt}}
\end{equation*}
\end{enumerate}

The second proof of Lemma~\ref{lem-column} is a bijection, which is
built on an involution  $\rho$  on the fillings of a rectangular
shape $\M$.

\noindent \textbf{An involution $\rho$ on rectangular shapes}. \\
Let $\M$ be an $n \times m$  rectangle.  We order the columns of
$\M$ from left to right, i.e., $C_1\prec' \cdots  \prec' C_m$, and
set $L(\M)=\M$. For any filling $M$, give it the coloring as
described in Subsection 4.1, and apply the bijection $\Psi$ from
$\F(\mathcal {M},\bve, \s)$ to $\mathcal{C}_{s_1+1}(h_1-s_1)\times
\cdots \times \mathcal{C}_{s_m+1}(h_m-s_m)$. For any filling $M$
with $\Psi(M)=(c^{(1)},c^{(2)},\ldots,c^{(m)})$ under the order
$\prec'$, let $\rho(M)$ be the filling whose associated sequence of
compositions is $(c^{(1),r},c^{(2)},\ldots,c^{(m)})$, again under
the order
 $\prec'$, where
\[
c^{(1),r}=(c_{s_1+1}^{(1)},\ldots,c_2^{(1)},c_1^{(1)}) \text{~if~}
c^{(1)}=(c_{1}^{(1)},c_2^{(1)},\ldots,c_{s_1+1}^{(1)} ). \] Then it
is easy to verify that $\rho(\rho(M))=M$ and $(\gamma(\{1\}; M),
\gamma(\overline{\{1\}};M))=(\se(\rho(M)), \ne(\rho(M)))$.

\noindent \emph{Second proof of Lemma \ref{lem-column}. }
 Given a general moon polyomino $\M$, assume that the rows intersecting the first column
are $\{R_a, \dots, R_b\}$. Let $\M_c$ be the union $R_a \cup \cdots
\cup R_b$. Clearly, for any $M \in \F(\M, \bve, \s)$, a left
$\{1\}$-NE (SE) chain  consists of two $1$-cells in $\M_c$.  Let
$C_i ' =C_i \cap \M_c$ be the restriction of the column $C_i$ on
$\M_c$. Then $C_1'=C_1$ and $|C_1'| \geq |C_2'| \geq \cdots \geq
|C_m'|$.

Suppose that
\begin{eqnarray*}
|C_1'|=|C_2'|=\cdots=|C_{j_1}'|>|C'_{j_1+1}|=|C'_{j_1+2}|= \cdots
=|C'_{j_2}| > |C'_{j_2+1}| \cdots \\[3pt]
 \cdots
=|C'_{j_{k-1}}|>|C'_{j_{k-1}+1}|=|C'_{j_{k-1}+2}|=\cdots
=|C'_{j_k}|=|C'_{m}|.
\end{eqnarray*}
Let $B_i$ be the greatest rectangle contained in $\M_c$
whose right most column is
        $C'_{j_{i}}$ ($1 \leq i \leq k$), and $B_i' = B_i \cap B_{i+1}$ ($1\leq i \leq k-1$).

We define $\phi_\gamma\colon \F(\M, \bve, \s)\rightarrow \F(\M,
\bve, \s)$ by constructing a sequence of fillings $(M, M_k, \dots,
M_1)$ starting from $M$.

\noindent \textbf{The map $\phi_\gamma\colon  \F(\M, \bve,
\s)\rightarrow \F(\M, \bve, \s)$}

Let $M \in \F(\M, \bve, \s)$.
\begin{enumerate}
 \item  The filling $M_k$ is obtained from $M$ by replacing $M \cap B_k$ with $\rho(M \cap B_k)$.
 \item  For $i$ from $k-1$ to $1$:
    \begin{enumerate}
       \item   Define a filling $N_i$ on $B_i'$ by setting $N_i =\rho(M_{i+1} \cap B_i')$.
 Let the filling $M_i'$ be obtained from $M_{i+1}$ by replacing $M_{i+1} \cap B_i'$ with  $N_i$.
       \item    The filling $M_i$ is obtained from $M_i'$ by replacing $M_i' \cap B_i$ with $\rho(M_i' \cap B_i)$.
      \end{enumerate}
\item Set $\phi_\gamma(M)=M_1$.
\end{enumerate}

See Example \ref{ex3} for an illustration.

\noindent {Claim}:  $(\gamma(\{1\};M), \gamma(\overline{\{1\}}; M))=
 (\se(\phi_\gamma(M)), \ne(\phi_\gamma(M)))$.

We are able to keep track of the statistic $\gamma(\{1\}; M)$ in the
above algorithm. In Step 1, by the definition of $\rho$ we have
\begin{eqnarray*}
 \gamma(\{1\}; M)& =& \#\{\text{left $\{1\}$-NE chain of $M$} \} + \#\{\text{left $\overline{\{1\}}$-SE chain of $M$}\} \\
                  & =& \#\{\text{left $\{1\}$-NE chain of $M_k$} \} + \#\{\text{left $\overline{\{1\}}$-SE chain of $M_k$}\} \\
                  &  & -\#\{\text{left $\{1\}$-NE chain of $M_k$ in $B_k$}\}
 + \#\{\text{left $\{1\}$-SE chain of $M_k$ in $B_k$}\}.
\end{eqnarray*}
Let $\mathcal{B}_i=B_i \cup \cdots \cup B_k$.
For $i$ from $k-1$ to $1$,  Step 2(a) implies that
for the filling $M_i'$,
\begin{eqnarray*}
 \gamma(\{1\}; M)& = &\#\{\text{left $\{1\}$-NE chain of $M_i'$} \} + \#\{\text{left $\overline{\{1\}}$-SE chain of $M_i'$}\} \\
                  &  & -\#\{\text{left $\{1\}$-NE chain of $M_i'$ in $\mathcal{B}_{i+1}$}\}
                   + \#\{\text{left $\{1\}$-SE chain of $M_i'$ in $\mathcal{B}_{i+1}$}\}  \\
                   &  & +\#\{\text{left $\{1\}$-NE chain of $M_i'$ in $B_i'$}\}
                        - \#\{\text{left $\{1\}$-SE chain of $M_i'$ in $B_i'$}\}.
\end{eqnarray*}
Then Step 2(b) implies that in the filling $M_i$,
\begin{eqnarray*}
 \gamma(\{1\}; M)  & =& \#\{\text{left $\{1\}$-NE chain of $M_i$} \} + \#\{\text{left $\overline{\{1\}}$-SE chain of $M_i$}\} \\
                  &  & -\#\{\text{left $\{1\}$-NE chain of $M_i$ in $\mathcal{B}_i$}\}
             +  \#\{\text{left $\{1\}$-SE chain of $M_i$ in $\mathcal{B}_i$}\}.
\end{eqnarray*}
Since all the $\{1\}$-NE (SE) chains of $M_i$ are in $\M_c=B_1\cup
\cdots \cup B_k=\mathcal{B}_1$, when $i=1$ we have $\gamma(\{1\}; M)
= \se(M_1)=\se(\phi_\gamma(M))$. Similarly,
$\gamma(\overline{\{1\}}; M)=\ne(\phi_\gamma(M))$. \qed

\begin{figure}[ht]
\begin{center}
\setlength{\unitlength}{0.35mm}
\begin{picture}(310,100)

\qbezier[500](20,20)(22,20)(30,20) \qbezier[500](0,30)(25,30)(50,30)
\qbezier[500](0,40)(30,40)(60,40) \qbezier[500](0,50)(30,50)(60,50)
\qbezier[500](0,60)(30,60)(60,60) \qbezier[500](0,70)(25,70)(50,70)
\qbezier[500](0,80)(15,80)(30,80) \qbezier[500](10,90)(15,90)(30,90)
\qbezier[500](0,30)(0,55)(0,80) \qbezier[500](10,30)(10,60)(10,90)
\qbezier[500](20,20)(20,55)(20,90)
\qbezier[500](30,20)(30,55)(30,90)
\qbezier[500](40,30)(40,50)(40,70)\qbezier[500](50,30)(50,50)(50,70)
\qbezier[500](60,40)(60,50)(60,60)

\linethickness{0.95pt}

 \put(0,70){\line(0,-1){40}} \put(50,70){\line(0,-1){40}}
\put(0,30){\line(1,0){50}}  \put(0,70){\line(1,0){50}}

\thinlines

\put(3,62){\footnotesize{1}} \put(13,82){\footnotesize{1}}
\put(23,22){\footnotesize{1}} \put(23,52){\footnotesize{1}}
\put(23,72){\footnotesize{1}} \put(33,32){\footnotesize{1}}
\put(53,42){\footnotesize{1}}

\put(3,10){\makebox(0,0)[l]{\footnotesize$M=M_3=M'_2$}}
\put(0,-2){\makebox(0,0)[l]{\footnotesize(boxed part is $B_2$)}}

\put(70,50){\vector(1,0){25}} \put(80,54){$\rho$}
\put(73,40){\footnotesize{to $B_2$}}

\put(130,20){\line(1,0){10}}   \put(110,30){\line(1,0){50}}
\put(110,40){\line(1,0){60}}   \put(110,50){\line(1,0){60}}
\put(110,60){\line(1,0){60}}   \put(110,70){\line(1,0){50}}
\put(110,80){\line(1,0){30}}   \put(120,90){\line(1,0){20}}
\put(110,80){\line(0,-1){50}}  \put(120,90){\line(0,-1){60}}
\put(130,90){\line(0,-1){70}}  \put(140,90){\line(0,-1){70}}
\put(150,70){\line(0,-1){40}}  \put(160,70){\line(0,-1){40}}
\put(170,60){\line(0,-1){20}}

\linethickness{0.95pt}

 \put(110,70){\line(0,-1){40}}
\put(140,70){\line(0,-1){40}} \put(110,30){\line(1,0){30}}
\put(110,70){\line(1,0){30}}

\thinlines

\put(113,32){\footnotesize{1}} \put(123,82){\footnotesize{1}}
\put(133,72){\footnotesize{1}} \put(133,62){\footnotesize{1}}
\put(133,22){\footnotesize{1}} \put(143,52){\footnotesize{1}}
\put(163,42){\footnotesize{1}}

\put(135,10){\makebox(0,0)[l]{\footnotesize$M_2$}}
\put(110,-2){\makebox(0,0)[l]{\footnotesize(boxed part is $B'_1$)}}

\put(180,50){\vector(1,0){25}} \put(190,54){$\rho$}
\put(183,40){\footnotesize{to $B'_1$}}
\put(240,20){\line(1,0){10}}   \put(220,30){\line(1,0){50}}
\put(220,40){\line(1,0){60}}   \put(220,50){\line(1,0){60}}
\put(220,60){\line(1,0){60}}   \put(220,70){\line(1,0){50}}
\put(220,80){\line(1,0){30}}   \put(230,90){\line(1,0){20}}
\put(220,80){\line(0,-1){50}}  \put(230,90){\line(0,-1){60}}
\put(240,90){\line(0,-1){70}}  \put(250,90){\line(0,-1){70}}
\put(260,70){\line(0,-1){40}}  \put(270,70){\line(0,-1){40}}
\put(280,60){\line(0,-1){20}}

\linethickness{0.95pt}

\put(220,80){\line(0,-1){50}} \put(250,80){\line(0,-1){50}}
\put(220,30){\line(1,0){30}}  \put(220,80){\line(1,0){30}}

\thinlines

\put(223,62){\footnotesize{1}} \put(233,82){\footnotesize{1}}
\put(243,72){\footnotesize{1}} \put(243,32){\footnotesize{1}}
\put(243,22){\footnotesize{1}} \put(253,52){\footnotesize{1}}
\put(273,42){\footnotesize{1}}

\put(235,10){\makebox(0,0)[l]{\footnotesize$M'_1=M_1$}}
\put(220,-2){\makebox(0,0)[l]{\footnotesize(boxed part is $B_1$)}}

\end{picture}
\end{center}

\caption{The map $\phi_\gamma$.} \label{phi_gamma}

\end{figure}
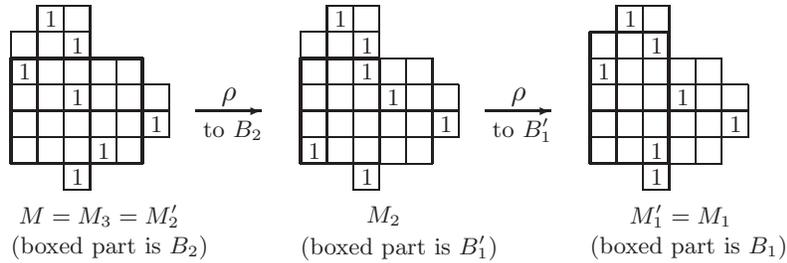

\begin{exa} \label{ex3}
 \emph{Figure \ref{phi_gamma} shows an example of the map $\phi_\gamma$ applied
to a filling $M$.  The filling $M$ is given in the figure on the
left, where $|C_1|=|C_2'|=|C_3'| > |C_4'| =|C_5'|>|C_6'|$. Hence
$k=3$, $j_1=3$, $j_2=5$ and $j_3=6$. It is easy to see that
$M=M_3=M_2'$. Figure \ref{phi_gamma} shows how to get $M_2$ and
$M_1'$. In this example, it happens that $M_1'=M_1$. }
\end{exa}

\begin{prop}\label{prop-column}
Assume  $T=\{c_1,c_2,\ldots, c_t\} \subseteq [m]$ with $c_1<c_2<\cdots<c_t$. Let
$T'=\{c_1,c_2,\ldots, c_{t-1}\}$. Then $G^l_T(p,q)=G^l_{T'}(p,q)$.
\end{prop}
\pf Like Lemma \ref{lem-column}, Prop.~\ref{prop-column} can be
proved either by analyzing the generating functions, or by a
bijection built on the map $\phi_\gamma$.
 Here we give the details  of the bijection which
will be used in Section 5.

Let $U=\{C_i: 1 \leq j < c_t\}$ be the set of columns on the left of
column $C_{c_t}$, and $V$ be the set of remaining columns. For any
$M \in \F(\M, \bve, \s)$, let $\xi_{c_t}(M)$ be the filling obtained
from $M$ by replacing $M \cap V$ with $\phi_\gamma(M \cap V)$.
Then $\xi_{c_t}$ is a bijection on
$\F(\M, \bve, \s)$ such that
\begin{eqnarray} \label{xi_t}
(\gamma(T;M),\gamma(\bar T;M)) =
(\gamma(T'; \xi_{c_t}(M)), \gamma(\bar{T'};\xi_{c_t}(M))).
\end{eqnarray}
The proof is similar to that of Prop.~\ref{prop-row} and is omitted.
\qed

\noindent \emph{Proof of Theorem \ref{thm-column}}. Assume $T=\{c_1,
c_2, \dots, c_t\} \subseteq \mathcal{C}$ with $c_1 < c_2 < \cdots <
c_t$. Let $\Sigma_\gamma= \xi_{c_1} \circ \xi_{c_2} \circ \cdots
\circ \xi_{c_t}$, where $\xi_c$ is  defined in the proof of Prop.
\ref{prop-column}. Then $\Sigma_\gamma$ is a bijection on $\F(\M,
\bve, \s)$  with the property that
\begin{equation*}
 (\gamma(T; M), \gamma(\bar T; M)) =
 (\se(\Sigma_\gamma(M)), \ne(\Sigma_\gamma(M))).  \tag*{\rule{4pt}{7pt}}
\end{equation*}

\section{Invariance Properties}

The bi-variate generating function of $(\ne, \se)$ (cf. Theorem
\ref{gf}) implies that the mixed statistics are invariant under any
permutation of rows and/or columns. To be more specific, let $\M$ be
a moon polyomino. For any moon polyomino $\M'$ obtained from $\M$ by
permuting the rows and/or the columns of $\M$, we have
\begin{eqnarray*}
&&\# \{ M \in \F(\M, \bve, \s): \lambda(A; M)=i, \lambda(\bar{A};
M)=j\}
\\[3pt]
&&\quad= \#\{ M' \in \F(\M', \bve', \s'): \lambda(A; M')=i,
\lambda(\bar{A}; M')=j\}
\end{eqnarray*}
for any nonnegative integers $i$ and $j$, where $\bve'$ (resp.
$\s'$) is the sequence obtained from $\bve$ (resp. $\s$) in the same
ways as the  rows (resp. columns) of $\M'$ are obtained from the
 rows (resp. columns) of $\M$, and $\lambda(A;M)$ is any of the
 four statistics $\alpha(S;M), \beta(S; M), \gamma(T;M)$, and $\delta(T;M)$.
In this section we present bijective proofs of such phenomena.

Let $\M$ be a general moon polyomino. Let $\N_l$ be the unique
left-aligned moon polyomino whose sequence of row lengths is equal
to $|R_1|, \ldots, |R_n|$ from top to bottom. In other words, $\N_l$
is the left-aligned polyomino obtained by rearranging the columns of
$\M$ by length in weakly decreasing order from left to right. We
shall use an algorithm developed in  \cite{Chen09} that rearranges
the columns of $\M$ to generate $\N_l$.

\noindent \textbf{The algorithm $\alpha$ for rearranging $\M$:}
\begin{itemize}
\item [Step 1] Set $\M'=\M$.
\item [Step 2] If $\M'$ is left aligned, go to Step 4.
\item [Step 3] If $\M'$ is not left-aligned, consider the largest rectangle
$\mathcal{B}$ completely contained in $\M'$ that contains $C_1$, the
leftmost column of $\M'$. Update $\M'$ by setting $\M'$ to be the
polyomino obtained by moving the leftmost column of $\mathcal{B}$ to
the right end. Go to Step 2.
\item [Step 4] Set $\N_l=\M'$.
\end{itemize}

\begin{figure}[ht]
\begin{center}
\setlength{\unitlength}{0.35mm}
\begin{picture}(400,70)

\put(20,0){\line(1,0){30}}  \put(0,10){\line(1,0){70}}
\put(0,20){\line(1,0){70}}  \put(0,30){\line(1,0){70}}
\put(0,40){\line(1,0){60}}  \put(10,50){\line(1,0){40}}
\put(20,60){\line(1,0){30}} \put(30,70){\line(1,0){20}}

\put(0,40){\line(0,-1){30}} \put(10,50){\line(0,-1){40}}
\put(20,60){\line(0,-1){60}}\put(30,70){\line(0,-1){70}}
\put(40,70){\line(0,-1){70}}\put(50,70){\line(0,-1){70}}
\put(60,40){\line(0,-1){30}}\put(70,30){\line(0,-1){20}}

\linethickness{0.95pt}

\put(0,40){\line(0,-1){30}} \put(60,40){\line(0,-1){30}}
\put(0,10){\line(1,0){60}}  \put(0,40){\line(1,0){60}}

\thinlines

\put(0,30){\line(1,1){10}} \put(0,35){\line(1,1){5}}
\put(5,30){\line(1,1){5}}

\put(0,20){\line(1,1){10}} \put(0,25){\line(1,1){5}}
\put(5,20){\line(1,1){5}}

\put(0,10){\line(1,1){10}} \put(0,15){\line(1,1){5}}
\put(5,10){\line(1,1){5}}

\put(80,35){\vector(1,0){20}}

\put(120,0){\line(1,0){30}}  \put(110,10){\line(1,0){70}}
\put(110,20){\line(1,0){70}} \put(110,30){\line(1,0){70}}
\put(110,40){\line(1,0){60}} \put(110,50){\line(1,0){40}}
\put(120,60){\line(1,0){30}} \put(130,70){\line(1,0){20}}

\put(110,50){\line(0,-1){40}} \put(120,60){\line(0,-1){60}}
\put(130,70){\line(0,-1){70}} \put(140,70){\line(0,-1){70}}
\put(150,70){\line(0,-1){70}} \put(160,40){\line(0,-1){30}}
\put(170,40){\line(0,-1){30}} \put(180,30){\line(0,-1){20}}

\linethickness{0.95pt}

\put(110,50){\line(0,-1){40}} \put(150,50){\line(0,-1){40}}
\put(110,10){\line(1,0){40}}  \put(110,50){\line(1,0){40}}

\thinlines

\put(110,40){\line(1,1){10}} \put(110,45){\line(1,1){5}}
\put(115,40){\line(1,1){5}}

\put(110,30){\line(1,1){10}} \put(110,35){\line(1,1){5}}
\put(115,30){\line(1,1){5}}

\put(110,20){\line(1,1){10}} \put(110,25){\line(1,1){5}}
\put(115,20){\line(1,1){5}}

\put(110,10){\line(1,1){10}} \put(110,15){\line(1,1){5}}
\put(115,10){\line(1,1){5}}

\put(160,40){\line(1,-1){10}} \put(160,35){\line(1,-1){5}}
\put(165,40){\line(1,-1){5}}

\put(160,30){\line(1,-1){10}} \put(160,25){\line(1,-1){5}}
\put(165,30){\line(1,-1){5}}

\put(160,20){\line(1,-1){10}} \put(160,15){\line(1,-1){5}}
\put(165,20){\line(1,-1){5}}

\put(160,30){\line(1,1){10}} \put(160,35){\line(1,1){5}}
\put(165,30){\line(1,1){5}}

\put(160,20){\line(1,1){10}} \put(160,25){\line(1,1){5}}
\put(165,20){\line(1,1){5}}

\put(160,10){\line(1,1){10}} \put(160,15){\line(1,1){5}}
\put(165,10){\line(1,1){5}}

\put(190,35){\vector(1,0){20}}

\put(220,0){\line(1,0){30}} \put(220,10){\line(1,0){70}}
\put(220,20){\line(1,0){70}}\put(220,30){\line(1,0){70}}
\put(220,40){\line(1,0){60}}\put(220,50){\line(1,0){40}}
\put(220,60){\line(1,0){30}}\put(230,70){\line(1,0){20}}

\put(220,60){\line(0,-1){60}}\put(230,70){\line(0,-1){70}}
\put(240,70){\line(0,-1){70}}\put(250,70){\line(0,-1){70}}
\put(260,50){\line(0,-1){40}}\put(270,40){\line(0,-1){30}}
\put(280,40){\line(0,-1){30}}\put(290,30){\line(0,-1){20}}

\linethickness{0.95pt}

\put(220,60){\line(0,-1){60}} \put(250,60){\line(0,-1){60}}
\put(220,0){\line(1,0){30}}  \put(220,60){\line(1,0){30}}

\thinlines

\put(220,50){\line(1,1){10}} \put(220,55){\line(1,1){5}}
\put(225,50){\line(1,1){5}}

\put(220,40){\line(1,1){10}} \put(220,45){\line(1,1){5}}
\put(225,40){\line(1,1){5}}

\put(220,30){\line(1,1){10}} \put(220,35){\line(1,1){6}}
\put(225,30){\line(1,1){5}}

\put(220,20){\line(1,1){10}} \put(220,25){\line(1,1){5}}
\put(225,20){\line(1,1){5}}

\put(220,10){\line(1,1){10}} \put(220,15){\line(1,1){5}}
\put(225,10){\line(1,1){5}}

\put(220,0){\line(1,1){10}} \put(220,5){\line(1,1){5}}
\put(225,0){\line(1,1){5}}

\put(250,50){\line(1,-1){10}} \put(250,45){\line(1,-1){5}}
\put(255,50){\line(1,-1){5}}

\put(250,40){\line(1,-1){10}} \put(250,35){\line(1,-1){5}}
\put(255,40){\line(1,-1){5}}

\put(250,30){\line(1,-1){10}} \put(250,25){\line(1,-1){5}}
\put(255,30){\line(1,-1){5}}

\put(250,20){\line(1,-1){10}} \put(250,15){\line(1,-1){5}}
\put(255,20){\line(1,-1){5}}

\put(250,40){\line(1,1){10}} \put(250,45){\line(1,1){5}}
\put(255,40){\line(1,1){5}}

\put(250,30){\line(1,1){10}} \put(250,35){\line(1,1){5}}
\put(255,30){\line(1,1){5}}

\put(250,20){\line(1,1){10}} \put(250,25){\line(1,1){5}}
\put(255,20){\line(1,1){5}}

\put(250,10){\line(1,1){10}} \put(250,15){\line(1,1){5}}
\put(255,10){\line(1,1){5}}

\put(300,35){\vector(1,0){20}}

\put(330,0){\line(1,0){30}} \put(330,10){\line(1,0){70}}
\put(330,20){\line(1,0){70}}\put(330,30){\line(1,0){70}}
\put(330,40){\line(1,0){60}}\put(330,50){\line(1,0){40}}
\put(330,60){\line(1,0){30}}\put(330,70){\line(1,0){20}}

\put(330,70){\line(0,-1){70}}\put(340,70){\line(0,-1){70}}
\put(350,70){\line(0,-1){70}}\put(360,60){\line(0,-1){60}}
\put(370,50){\line(0,-1){40}}\put(380,40){\line(0,-1){30}}
\put(390,40){\line(0,-1){30}}\put(400,30){\line(0,-1){20}}

\put(350,60){\line(1,-1){10}} \put(350,55){\line(1,-1){5}}
\put(355,60){\line(1,-1){5}}

\put(350,50){\line(1,-1){10}} \put(350,45){\line(1,-1){5}}
\put(355,50){\line(1,-1){5}}

\put(350,40){\line(1,-1){10}} \put(350,35){\line(1,-1){5}}
\put(355,40){\line(1,-1){5}}

\put(350,30){\line(1,-1){10}} \put(350,25){\line(1,-1){5}}
\put(355,30){\line(1,-1){5}}

\put(350,20){\line(1,-1){10}} \put(350,15){\line(1,-1){5}}
\put(355,20){\line(1,-1){5}}

\put(350,10){\line(1,-1){10}} \put(350,5){\line(1,-1){5}}
\put(355,10){\line(1,-1){5}}

\put(350,50){\line(1,1){10}} \put(350,55){\line(1,1){5}}
\put(355,50){\line(1,1){5}}

\put(350,40){\line(1,1){10}} \put(350,45){\line(1,1){5}}
\put(355,40){\line(1,1){5}}

\put(350,30){\line(1,1){10}} \put(350,35){\line(1,1){5}}
\put(355,30){\line(1,1){5}}

\put(350,20){\line(1,1){10}} \put(350,25){\line(1,1){5}}
\put(355,20){\line(1,1){5}}

\put(350,10){\line(1,1){10}} \put(350,15){\line(1,1){5}}
\put(355,10){\line(1,1){5}}

\put(350,0){\line(1,1){10}} \put(350,5){\line(1,1){5}}
\put(355,0){\line(1,1){5}}

\end{picture}
\end{center}
\caption{The algorithm $\alpha$.} \label{alpha}
\end{figure}
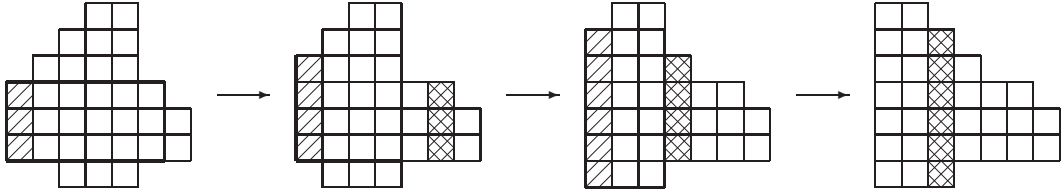

Figure~\ref{alpha} is an illustration of the algorithm $\alpha$.

Based on the algorithm $\alpha$, Chen et al. constructed a bijection
$g=g_{\M}: \F(\M, \bve, \s)\rightarrow \F(\N_l, \bve, \s')$ such
that $(\se(M), \ne(M))  =(\se(g(M)), \ne(g(M)))$, see \cite[Section
5.3.2]{Chen09}.

Combining $g_{\M}$ with the bijection $\Theta_\alpha$ constructed in
the proof of Theorem \ref{thm-row}, we are led to the following
invariance property.

\begin{thm}
Let $\M$ be a moon polyomino. For any moon polyomino $\M'$ obtained from
$\M$ by permuting the columns of $\M$, the map
\begin{eqnarray} \label{comp1}
\Phi_\alpha = \Theta_\alpha^{-1} \circ g_{\M'}^{-1} \circ g_{\M} \circ \Theta_\alpha:
\F(\M, \bve, \s) \rightarrow \F(\M', \bve, \s')
\end{eqnarray}
is a bijection with the property that
\[
(\alpha(S; M), \alpha(\bar S;M) )= (\alpha(S; M'), \alpha(\bar S;M')).
\]
\end{thm}

Similarly, let $\N_t$ be the top aligned polyomino obtained from
$\M$ by rotating 90 degrees counterclockwise first, followed by
applying the algorithm $\alpha$, and finally rotating 90 degrees
clockwise. Such operations enable us to establish a bijection
$h=h_{\M}$ from $\F(\M, \bve, \s)$ to $\F(\N_t, \bve', \s)$ that
keeps the statistics $(\se, \ne)$.  The map $h_{\M}$ can be
described by using the map $g_{\M}$ under the algorithm $\alpha$
with a rotation of 90 degrees clockwise. More precisely, the rotated
algorithm $\alpha'$ is the same as the algorithm $\alpha$, except
that the term \emph{left-aligned} is replaced with the term
\emph{top-aligned}, $C_1$ is replaced with $R_1$, and \emph{left}
and \emph{right} are replaced with \emph{top} and \emph{bottom}
respectively. In fact the map $h_{\M}$ is much simpler than $g_{\M}$
since every row in the filling has at most one $1$-cell. We state it
in full  detail for  completeness.

\begin{thm}
There is a bijection $h_{\M}\colon \F(\M, \bve, \s)\rightarrow
\F(\N_t, \bve', \s)$ such that $(\se(M), \ne(M))=(\se(h(M)),
\ne(h(M)))$.
\end{thm}

\pf Let $M \in \F(\M, \bve, \s)$. To obtain $h_{\M}(M)$, we perform
the rotated algorithm $\alpha'$ to transform the shape $\M$ to
$\N_t$ and change the filling when we move rows down in Step 3 so
that the number of 1's in each row and column is preserved.

Let $N$ be the filling on the rectangular $\B$ in Step 3 of the
rotated  algorithm $\alpha'$ that contains the row $R_1$ of the
current filling.  Let $\B'$ be the rectangle obtained by moving the
row $R_1$ from top to the bottom  of $\B$. Fill it to obtain a
filling $N'$ as follows.

 1. If $R_1$ is empty, then $N'$ is obtained from $N \setminus \{R_1\}$
by adding an empty row at the bottom.

 2. If $R_1$ has a $1$-cell,  \\
\mbox{}\hspace{1cm} (a) The rows that are empty in $\B$ remain empty
in $\B'$.
Shade these rows in both $\B$ and $\B'$.   \\
\mbox{}\hspace{1cm} (b) The filling on the rectangle formed by the
un-shaded rows of $\B'$ is the same as
         $N$ restricted to the rectangle obtained from the un-shaded rows of $\B$.

The filling outside $\B$ remains unchanged.

Applying the rotated algorithm $\alpha'$  with the above operations
on filling $M$, we finally obtain the filling $h_{\M}(M)$. The proof
of \cite[Prop. 5.10]{Chen09} ensures that $h_{\M}$ is a
 bijection. \qed

Combining the bijection $\Theta_\alpha$ with $h_{\M}$, we arrive at
the second invariance property.

\begin{thm}
Let $\M$ be a moon polyomino. For any moon polyomino $\M'$ obtained from
$\M$ by permuting the rows of $\M$, the map
\begin{eqnarray}\label{comp2}
\Lambda_\alpha = \Theta_\alpha^{-1} \circ h_{\M'}^{-1} \circ h_{\M} \circ \Theta_\alpha:
\F(\M, \bve, \s) \rightarrow \F(\M', \bve', \s)
\end{eqnarray}
is a bijection with the property that
\[
(\alpha(S; M), \alpha(\bar S;M) )= (\alpha(S; M'), \alpha(\bar
S;M')).
\]
\end{thm}

It is evident that replacing $\Theta_\alpha$ with the map
$\Sigma_\gamma$ defined in the proof of Theorem \ref{thm-column} in
 \eqref{comp1} and \eqref{comp2} leads to bijections preserving
the statistics $(\gamma(T; M), \gamma(\bar T; M))$ under any
permutation of columns or rows. Similar results hold for the
statistics $\beta(S; M)$ and $\delta(T; M)$ by reflecting  the moon
polyomino with respect to a horizonal or a vertical line.


\end{document}